\documentclass[12pt]{article}
\usepackage{a4wide}
\usepackage{amssymb}
\usepackage{times}

\begin{document}
\newtheorem{theorem}{Theorem}[section]
\newtheorem{definition}[theorem]{Definition}
\newtheorem{lemma}[theorem]{Lemma}
\newtheorem{proposition}[theorem]{Proposition}
\newtheorem{remark}[theorem]{Remark}
\newtheorem{remarks}[theorem]{Remarks}
\newtheorem{examples}[theorem]{Examples}
\newtheorem{example}[theorem]{Example}
\newtheorem{corollary}[theorem]{Corollary}
\newtheorem{conjecture}[theorem]{Conjecture}

\newcommand{\IC}{\mathbb{C}}
\newcommand{\IR}{\mathbb{R}}
\newcommand{\IZ}{\mathbb{Z}}
\def\qedrule{\vrule height0.75em width0.5em depth0.25em}
\def\makeqedrule{\nobreak \ifvmode \relax \else
      \ifdim\lastskip<1.5em \hskip-\lastskip
		  \hskip1.5em plus0em minus0.5em \fi \nobreak
	\qedrule\fi}
\let\qed=\makeqedrule

\def\blb#1{\mathbb{#1}}
\def\Bbb#1{\mathbb{#1}}

\title{Complex Geodesics on Convex Domains          }
\author{Se\'{a}n Dineen and Richard M. Timoney         }        
\date{\normalsize
Dedicated to M. Valdivia on the occasion of his 60$^{\rm th}$
birthday
}

\vspace*{-23mm}
This paper appeared in Progress in Functional Analysis, K.D. Bierstedt, J.
Bonet, J. Horv\'ath \& M. Maestre (Eds.) pp 333--365, Elsevier (1992) .

\pagestyle{myheadings}
\markboth{S. Dineen, R.M. Timoney}{Comlex geodesics on complex
domains}
\thispagestyle{plain}
\vspace*{1.5cm}
\noindent
{\Large Complex Geodesics on Convex Domains          }\\[6mm]
{\large Se\'{a}n Dineen$^{\rm a}$ and Richard M. Timoney$^{\rm b}$}\\[6mm]
{\small Dedicated to M. Valdivia on the occasion of his 60$^{\rm th}$
birthday}\\[6mm]
{$^{\rm a}$\large Department of Mathematics,
University College Dublin,
Belfield,
Dublin 4,
Ireland}\\[6mm]
{$^{\rm b}$\large School of Mathematics,
39 Trinity College,
Dublin 2,
Ireland}
\vspace{6mm}

\section*{Abstract}
Existence and uniqueness of complex geodesics joining two
points of a convex bounded domain in a Banach space $X$ are considered.
Existence is proved for the unit ball of $X$ under the assumption that $X$
is 1-complemented in its double dual. Another existence result for taut
domains is also proved. Uniqueness is proved for strictly
convex bounded domains in spaces with the analytic Radon-Nikodym property.
If the unit ball of $X$ has a modulus of complex uniform convexity with
power type decay at 0, then all complex geodesics in the unit ball
satisfy a Lipschitz condition. The results are applied to classical Banach
spaces and to give a formula
describing all complex geodesics in the unit ball of the sequence spaces
$\ell^p$ ($1 \leq p < \infty$).
\vspace{1cm}

In this article, we discuss the existence, uniqueness and continuity
of complex
geodesics on a convex domain $\cal D$ in a complex Banach space $X$.
The term `complex geodesic' is due to Vesentini~\cite{vesentini},
although the concept was discussed by Carath\'{e}odory~\cite{carath} and
Reiffen~\cite{reiffen} under the name `metric plane'. Recent results
on this topic are to be found in
\cite{DTV,gent-nu,gent-bal,gent-l1,vesent-hol,vesent-inv,vigue-geod,vigue-q}.
Applications of complex geodesics to the study of biholomorphic
automorphisms and to fixed point sets are to be found in
\cite{carath,vesentini,vesent-hol,vesent-inv,vigue-geod,vigue-q}.

Our results on the existence
problem depend on topological properties of the Banach space $X$,
the results on uniqueness depend on the geometry of the boundary
$\partial \cal D$ and on an analytic-geometric property of $X$ (the
analytic Radon-Nikodym property), while the continuity ({\em i.e.}
continuous extensions to the boundary) is obtained using complex uniform
convexity.

In section~\ref{sec1}, we introduce complex geodesics and related
concepts and prove some basic results. In section~\ref{sec-exth},
we show that every pair of points in the unit ball $B_X$ of $X$
can be joined by a complex geodesic provided $X$ is 1-complemented
in its double dual $X^{**}$. In proving this we show that the
canonical embedding of $B_X$ in $B_{X^{**}}$ is an isometry
for the Kobayashi metrics. We also obtain  a new simple proof of
a result of Davie and Gamelin~\cite{DG} that bounded analytic
functions on $B_X$ extend to $B_{X^{**}}$.
In section~\ref{sec-uniq}, we use the analytic Radon-Nikodym property
and extreme points to obtain uniqueness results and in
Section~\ref{sec-cont} we prove continuity properties of complex geodesics.
In Section~\ref{sec-ex} we apply the results of the preceding three 
sections to $\ell^p$ and related spaces.

An examination of the methods of section~\ref{sec-exth} shows that
topological considerations are only used to obtain `Montel type' theorems
for mappings defined on the unit disc in $\IC$ with values in
the domain ${\cal D} \subset X$ under consideration. This led us
to consider (in an infinite-dimensional setting) the concepts of
taut and complete hyperbolic domains (section~\ref{sec-taut}). In
section~\ref{sec-taut}, we also apply our results to show that certain domains
have constant negative curvature in the Kobayashi metric.

A number of the results in this paper were announced in \cite{dineen}.
We refer to \cite{dineen} for background results and further details
on complex geodesics.

\section{COMPLEX GEODESICS}
\label{sec1}

$\Bbb D$ is the open unit disc in $\IC$ and $\rho$ will denote the Poincar\'{e}
distance on $\Bbb D$, i.e.
\[
\rho (z,w) = \tanh^{-1} \left( \left| \frac{z-w}{1 - \bar{w}z} \right|
\right) \qquad (z,w \in {\Bbb D}).
\]
The infinitesimal Poincar\'{e} distance $\alpha$ is defined by $\alpha(z,v) =
|v|/(1 - |z|^2)$ (for $v \in \IC$ and $z \in \Bbb D$).

For ${\cal D}_1$ and ${\cal D}_2$ domains in complex Banach spaces,
$H({\cal D}_1, {\cal D}_2 ) $ denotes the set of all ${\cal D}_2$-valued
holomorphic functions on ${\cal D}_1$. For a domain $\cal D$ in a
complex Banach space $X$, $p,q \in {\cal D}$ and $v \in V$,
\begin{eqnarray*}
C_{\cal D} (p,q) &=& \sup \{ \rho (f(p), f(q)) : f \in H({\cal D}, {\Bbb D})
\}\\
c_{\cal D} (p,v) &=& \sup \{ |f'(p)(v)|) : f \in H({\cal D}, {\Bbb D})
\}\\
\delta_{\cal D} (p,q) &=& \inf \{ \rho(u,v) : \exists f \in H( {\Bbb D},
{\cal D}), f(u) = p, f(v) = q \}\\
K_{\cal D} (p,q) &=& \inf \left\{ \sum_{i=1}^n \delta_{\cal D} (w_{i-1}, w_i) :
n \geq 1, p=w_0, w_1, \ldots, w_n = q \in {\cal D} \right\}\\
k_{\cal D} (p, v) &=& \inf \{ \eta > 0 : \exists f \in H({\Bbb D}, {\cal D}),
f(0) = p, f'(0)\eta = v \}.
\end{eqnarray*}
$C_{\cal D}$ is called the Carath\'{e}odory distance on $\cal D$, $K_{\cal
D}$ the Kobayashi distance and $c_{\cal D}$ and $k_{\cal D}$ are the
corresponding infinitesimal metrics.

For the unit disc, $C_{\Bbb D} = K_{\Bbb D} = \rho$ and $c_{\Bbb D} =
k_{\Bbb D} = \alpha$. In general $C_{\cal D} \leq K_{\cal D} \leq
\delta_{\cal D}$. Unlike $C_{\cal D}$ and $K_{\cal D}$,
$\delta_{\cal D}$ does not obey the
triangle inequality in general.
$K_{\cal D}$ is the largest distance function on
$\cal D$ smaller than $\delta_{\cal D}$.
Holomorphic mappings are contractions
relative to any one of the above distances or infinitesimal metrics.

\begin{definition}
\label{geod-def}
Let ${\cal D}$ be a domain in a complex Banach space $X$ and let $d$ be
a distance on $\cal D$. A mapping $\phi \in H( {\Bbb D} , {\cal D})$
is called a {\em complex $d$-geodesic} if
\[
\rho (u,v) = d(\phi(u), \phi(v) ) \mbox{ for } u,v \in {\Bbb D}.
\]
If $z,w \in \phi({\Bbb D})$ are distinct points, then we refer to $\phi$
as a complex $d$-geodesic joining $z$ and $w$.

We will use the term {\em complex geodesic} for `complex
$C_{\cal D}$-geodesic'.
\end{definition}

\begin{proposition}
\label{geod-equiv}
For a domain $\cal D$ in a Banach space and $\phi \in H({\Bbb D}, {\cal D})$,
the following are equivalent
\begin{description}
\item{(a)} $\phi$ is a complex geodesic.
\item{(b)} there exist distinct points $u,v \in \Bbb D$ such that
$\rho(u,v) = C_{\cal D} ( \phi(u), \phi(v) )$.
\item{(c)} there exists a point $z \in \Bbb D$ such that $\alpha(z,1) =
c_{\cal D} ( \phi(z), \phi'(z) )$.
\item{(d)} $\phi$ is biholomorphic from $\Bbb D$ to
an analytic set $\phi( {\Bbb D} )$ and $\phi( {\Bbb D} )$ is a holomorphic
retract of $\cal D$ (i.e. there exists $f \in H( {\cal D}, {\cal D})$
such that $f \circ f = f$ and $f({\cal D}) = \phi({\Bbb D})$).
\item{(e)} $\phi$ is a complex $K_{\cal D}$-geodesic and
\[
\delta_{\cal D} \big| \phi({\Bbb D}) =
K_{\cal D} \big| \phi({\Bbb D}) =
C_{\cal D} \big| \phi({\Bbb D}).
\]
\end{description}
\end{proposition}

\noindent
{\sc Proof}: The equivalence of (a), (b) and (c) is due to
Vesentini~\cite{vesent-hol} and the fact that (d) is equivalent to (a)
is due to Reiffen~\cite[p. 19]{reiffen}
(see also
Lempert~\cite{lempert-hol}).
(e) $ \Rightarrow$ (a) and it remains to show that (a) $\Rightarrow$ (e).

Fix $u \neq v \in \Bbb D$ and write $p = \phi(u)$, $q = \phi(v)$.
Since $\rho(u,v) = C_{\cal D}(p,q)$, Montel's theorem implies the
existence of $f \in H( {\cal D}, {\Bbb D})$ such that
$f(p) = u$ and $f(q) = v$. Now $f \circ \phi \in H( {\Bbb D}, {\Bbb D})$
and
\[
\rho ( f \circ \phi (u), f \circ \phi(v) ) = C_{\cal D} (p,q) = \rho (u,v).
\]
By the Schwarz-Pick lemma (see \cite[p.5]{dineen}), $f \circ \phi$ is
a biholomorphic automorphism of $\Bbb D$ and
$ \rho( f \circ \phi (z), f \circ \phi (w) ) = \rho (z,w)$ for all
$z,w \in \Bbb D$.
Since holomorphic mappings are contractions, we have
\begin{eqnarray}
\label{star}
\rho(z, w) &=& \rho ( f \circ \phi (z), f \circ \phi (w) ) \leq
                  C_{\cal D} ( \phi(z), \phi(w) ) \nonumber\\
           &\leq& K_{\cal D} ( \phi(z), \phi(w) ) \leq
                  \delta_{\cal D} ( \phi(z), \phi(w) ) \\
           &\leq& \rho (z,w) \nonumber
\end{eqnarray}
and (\ref{star}) consists entirely of equalities. Thus (a) implies (e). \qed

\begin{definition}
\label{C-connect-def}
We call a domain $\cal D$ in a complex Banach space  {\em C-connected}
if every pair of points of $\cal D$ can be joined by a complex
geodesic.
\end{definition}

We remark that it follows from results of Vigu\'{e}~\cite{vigue-q}
(see also \cite[proposition 11.15 and corollary 11.17]{dineen})
that a subset of a bounded convex finite dimensional domain $\cal D$
is the range of a complex geodesic if and only if it is a connected
one-dimensional analytic subset and a holomorphic retract of $\cal D$.

An immediate consequence of the equivalence of (a) and (e) in
Proposition~\ref{geod-equiv} is the following result.

\begin{proposition}
\label{C=K=delta}
If $\cal D$ is a C-connected domain in a Banach space, then
\[
C_{\cal D} = K_{\cal D} = \delta_{\cal D}.
\]
\end{proposition}

This leads to many examples of domains which are {\em not}
C-connected (see for instance \cite[p.103]{dineen}). In particular
a proper domain in $\IC$ is C-connected if and only if it is
simply connected.

On the other hand Lempert~\cite{lempert-bull,lempert-hol} (see also
Royden-Wong~\cite{Roy-W}) has shown that $C_{\cal D} = K_{\cal D}
= \delta_{\cal D}$ and that $c_{\cal D} = k_{\cal D}$ if $\cal D$
is a convex bounded domain in $\IC^n$. This result was extended to
convex domains in arbitrary Banach spaces in \cite{DTV}, where it was
used to prove that the following are C-connected (see also
\cite[pp.90--91]{dineen})
\begin{description}
\item{(a)} bounded convex domains in reflexive Banach spaces
\item{(b)} the open unit ball $B_X$ of a dual Banach space $X$.
\end{description}

These facts lead us to propose the following conjecture.
\begin{conjecture}
\label{conject}
If $\cal D$ is a domain in a complex Banach space $X$ which
is biholomorphically equivalent to a bounded domain,
then the following are equivalent:
\begin{description}
\item{(a)} $C_{\cal D} = K_{\cal D} = \delta_{\cal D}$
\item{(b)} $\cal D$ is biholomorphically equivalent to a convex domain
\item{(c)} $\cal D$ is C-connected.
\end{description}
\end{conjecture}

(Note that (c) implies (a) and (b) implies (a).)

\section{EXTENSION AND EXISTENCE THEOREMS}
\label{sec-exth}
\setcounter{equation}{0}
In this section,
we consider a complex Banach space $X$ as a subspace of its double dual
space $X^{**}$ via the natural embedding. $B_X$ denotes the open unit ball
of $X$ and $H^\infty (B_X)$ denotes the space of scalar-valued bounded
holomorphic functions on $B_X$, with the supremum norm. $H^\infty(B_X)$
is a Banach algebra. We abbreviate $K_{B_X}$ as $K_X$ and $k_{B_X}$ as
$k_X$ from now on.

If $\cal U$ is an ultrafilter on a set $I$,
then $(X)_{\cal U}$ will denote the ultrapower of a Banach space $X$.
More specifically, if $\ell^\infty (I, X)$ denotes the space of bounded
$X$-valued functions on $I$ (with the supremum norm) and
\[
N_{\cal U} = \{ (x_i)_{i \in I} \in \ell^\infty (I, X) :
\lim_{\cal U} \| x_i \| = 0 \},
\]
then $(X)_{\cal U}$ is the quotient space $\ell^\infty (I, X) /
N_{\cal U}$.

There is a canonical isometric embedding $j \colon X \to (X)_{\cal U}$
given by
\[
j(x) = (x)_{i \in I} + N_{\cal U}
\]
(where $(x)_{i \in I}$ denotes the constant function $x$).

 The principle of local reflexivity in its ultrapower
formulation (see \cite{heinrich}) asserts
that given a Banach space $X$, there exists an ultrapower $(X)_{\cal U}$
such that
\begin{description}
\item{(i)}
there is an isometric embedding $J \colon X^{**} \to (X)_{\cal U}$
which extends the canonical embedding $j \colon X \to (X)_{\cal U}$;
\item{(ii)} the map
\begin{eqnarray*}
Q \colon (X)_{\cal U} &\to& X^{**}\\
Q( (x_i)_{i \in I} + N_{\cal U} ) &=& w^* \mbox{-} \lim_{\cal U} x_i
\end{eqnarray*}
satisfies $QJ = \mbox{id}_{X^{**}}$, $\|Q \| = 1$.

(Thus $JQ$ is a contractive projection of $(X)_{\cal U}$ onto the
isometric copy $J(X^{**})$ of $X^{**}$.)
\end{description}

The next lemma provides a crude version of the Schwarz inequality,
which we prove for completeness.

\begin{lemma}
\label{schwarz-crude}
Suppose $X$ is a complex Banach space, $f \in H(B_X, {\Bbb D})$,
$0 < \varepsilon < 1$ and $a, b \in (1- \varepsilon) B_X$. Then there
is a constant $A_\varepsilon$ depending only on $\varepsilon$ such
that
\[
\rho ( f(a), f(b) ) \leq A_\varepsilon \| a - b\|.
\]
\end{lemma}

\noindent
{\sc Proof}: We can first dispose of the case where $a$ and $b$ are
far apart, specifically the case $\| a - b \| > \varepsilon$. Define
$\phi \colon {\Bbb D} \to B_X$ by
\[
\phi(z) = \frac{z}{1-\varepsilon} a.
\]
Then $f \circ \phi \in H({\Bbb D}, {\Bbb D})$ and it follows
from the classical Schwarz-Pick lemma that
\[
\rho (f(a), f(0)) = \rho( f\circ \phi( 1-\varepsilon), f \circ \phi(0) )
\leq \rho( 0, 1 -\varepsilon).
\]
Since the same estimate also applies to $b$
\[
\rho ( f(a) , f(b) ) \leq 2 \rho( 0, 1 -\varepsilon) \leq A_\varepsilon
\varepsilon \leq A_\varepsilon \| a - b \|
\]
for $A_\varepsilon$ chosen suitably large (depending only on $\varepsilon$).

For $\| a - b \| \leq \varepsilon$, we can define $\phi \in H({\Bbb D},
B_X)$ by
\[
\phi(z) = \frac{a + b}{2} + \frac{\varepsilon (a-b)}{\|a-b\|} z
\]
and we will then have
\[
\phi \left( \frac{ \|a-b\|}{2\varepsilon} \right) = a, \quad
\phi \left( - \frac{ \|a-b\|}{2\varepsilon} \right) = b.
\]
Applying the Schwarz lemma to $f\circ \phi$, we find
that
\[
\rho( f(a) , f(b) ) \leq \rho \left( \|a-b\|/2\varepsilon, -\|a-b\|/2\varepsilon
\right)=
2 \rho ( \|a-b\|/2\varepsilon, 0) \leq A_\varepsilon \|a-b\|
\]
since $\|a-b\|/2\varepsilon \leq 1/2$. \qed

\begin{theorem}
\label{extend-ultra}
Suppose $Y = (X)_{\cal U}$ is an ultrapower of a complex Banach space
$X$, where $\cal U$ is an ultrafilter on the set $I$.

\begin{description}
\item{(i)} Suppose, for each $i\in I$,
 $f_i \colon B_X \to D$ is a holomorphic function from $B_X$
to the unit disc $D$ in the complex plane, and
\begin{equation}
\label{hypoth-limit}
\lim_{\cal U} | f_i (a_i) | < 1
\end{equation}
for one point $(a_i)_{i \in I} + N_{\cal U}$ of $B_Y$. Then the function
\begin{eqnarray*}
F \colon B_Y & \to & D \\
F( (x_i)_{i \in I} + N_{\cal U} ) &=& \lim_{\cal U} f_i (x_i)
\end{eqnarray*}
is well-defined and holomorphic.

\item{(ii)} The Kobayashi distance $K_Y$ on $B_Y$ is given by
\[
K_Y ( (x_i)_{i \in I}  + N_{\cal U}, (y_i)_{i \in I}  + N_{\cal U} )
= \lim_{\cal U} K_X( x_i, y_i).
\]
\end{description}
\end{theorem}

\noindent
{\sc Proof}: (i) Note first of all that if
$(x_i)_{i \in I}  + N_{\cal U}  \in B_Y$, then
$\lim_{\cal U} \|x_i\| < 1$ and so there exists
$U \in {\cal U}$ so that $\|x_i\| < 1$ for all $i \in U$.
Thus $f_i(x_i)$ is defined for $i \in U$ and $\lim_{\cal U} f_i(x_i)$
makes sense. In fact, it is possible to change $x_i$ for $i \not\in U$
so as to ensure $\sup_{i \in I} \|x_i\| < 1$
without changing the coset $(x_i)_{i \in I} + N_{\cal U}$
or the value of the limit.
We will make this change without comment
from now on for all points in $B_Y$ and,
in particular we will assume that $S = \sup_i \|a_i\| < 1$.

By compactness, $\lim_{\cal U} f_i(x_i)$ certainly exists in the closed
unit disc $\bar{\Bbb D}$. To show it is in the open disc fix $(x_i)_i$
and choose $0 < \varepsilon < 1-S$ so that $\sup_i \| x_i\| < 1- \varepsilon$.
By Lemma~\ref{schwarz-crude},
\[
\rho( f_i(x_i), f_i(a_i) ) \leq 2 A_\varepsilon \quad (i \in I).
\]
Choose now $U \in {\cal U}$ so that $T = \sup_{i \in U} |f_i(a_i)| < 1$.
It follows by the triangle inequality that
\[
\rho (f_i(x_i), 0) \leq 2 A_\varepsilon + \rho(0,f_i(a_i))
\leq 2 A_\varepsilon + \rho (0,T).
\]
Hence $\sup_{i \in U} |f_i(x_i)| < 1$ and $F$ does indeed map $B_Y$ into
$\Bbb D$.

Next we check that $F$ is continuous. For this, fix $(x_i)_i + N_{\cal U}
\in B_Y$ and choose $\varepsilon >0$ with $x = \sup_i \|x_i\| < 1 - 2 \varepsilon$.
For $y = (y_i)_i + N_{\cal U} \in B_Y$ and $\|x-y\|_Y < \varepsilon$,
we have $\|y_i\|_X < 1 - \varepsilon$ and so we
can apply Lemma~\ref{schwarz-crude} to see that
\[
\rho (f_i (x_i), f(y_i) ) \leq A_\varepsilon \| x_i - y_i\|
\]
Taking limits along $\cal U$, we deduce that
\[
\rho (F(x), F(y) ) \leq A_\varepsilon \| x - y\|_Y,
\]
which is enough to show continuity of $F$ at $x$.

Finally analyticity of $F$ follows from continuity together with
analyticity of the restricton of $F$ on complex lines in $Y$.
For $x = (x_i)_i + N_{\cal U} \in B_Y$ and $y = (y_i)_i +N_{\cal U} \in Y$,
Montel's theorem shows that
\[
F( x + zy) = \lim_{\cal U} f_i(x_i + zy_i)
\]
is an analytic function of $z$ on $\{ z \in \IC: \|x+zy\|_Y < 1 \}$.

\vspace{3mm}

\noindent
(ii) Fix $x = (x_i)_i + N_{\cal U}, y = (y_i)_i + N_{\cal U} \in B_Y$
and let $K_{\mbox{lim}} = \lim_{\cal U} K_X (x_i,y_i)$.
Choose $\varepsilon > 0$ and
put $r= \tanh (K_{\mbox{lim}}+ \varepsilon)$. Then $\rho(0,r) = K_{\mbox{lim}} + \varepsilon$. There
exists $U \in {\cal U}$ so that $K_X ( x_i,y_i) < K_{\mbox{lim}} + \varepsilon$ for
$i \in U$. Thus, since $\delta_X(x_i, y_i) = K_X (x_i, y_i)$, there
exists $g_i \in H ({\Bbb D}, B_X)$ satisfying $g_i(0) = x_i$,
$g_i(r) = y_i$ (for $i \in U$). For $i \in I \setminus U$, set $g_i =0$.

Now define $g \colon {\Bbb D} \to B_Y$, by $g(z) = ( g_i(z) )_i + N_{\cal U}$.
Of course we must check first that $\sup_i \| g_i(z) \|_X < 1$ for
all $z\in \Bbb D$. To this end, observe that
\begin{eqnarray*}
K_X ( g_i(z), 0) &\leq & K_X (g_i (z), g_i(0) ) + K_X(g_i(0), 0) \\
&\leq& \rho(z,0) + K_X(x_i,0)\\
&=& \rho(z,0) + \rho(\|x_i\|, 0)\\
&\leq & \rho(z,0) + \sup_{j \in I} \rho(\|x_j\|,0)\\
&< & \infty
\end{eqnarray*}
Here we have used the fact that holomorphic mappings are contractions
with respect to the Kobayashi distance and our standing assumption that
the $x_j$ are chosen so that $\sup_j \|x_j\| < 1$. The equality
$K_X(x,0) = \rho(\|x\|,0)$ (for $x \in B_X$) is elementary.
Hence we have
\[
\sup_i \rho( \|g_i(z)\|, 0) < \infty
\]
which implies $\sup_i \| g_i(z) \| < 1$.

Cauchy's formula shows that the functions $\|g_i''(z)\|$ are uniformly
bounded on compact subsets of $\Bbb D$. This will enable us to show by
a direct argument that $g$ is analytic. Fix $z \in \Bbb D$ and
suppose $|w-z| < (1-|z|)/2$. Then
\begin{eqnarray*}
\left\| \frac{g_i(w) - g_i(z)}{ w-z} - g_i'(z) \right\|
&=& \left\| \frac{1}{w-z} \int_z^w g'_i(\zeta) - g'_i(z) \, d\zeta\right\|\\
&=& \left\| \frac{1}{w-z} \int_z^w \int_z^\zeta g''_i(\eta) \, d\eta
           \, d\zeta\right\|\\
&\leq&  |w-z| \sup \{ |g''_i(\zeta)| : |\zeta| \leq (1+|z|)/2, i \in I\}
\end{eqnarray*}
Hence, if $\ell = (g'_i(z))_i + N_{\cal U}$, we have
\[
\left\| \frac{g(w) - g(z)}{w-z} - \ell \right\|_Y \leq M |w-z|
\]
for $|w-z| < (1-|z|)/2$ and $M$ a constant depending on $z$.
Taking the limit as $w \to z$, we see that $g'(z)$ exists and is $\ell$.

Therefore $g \in H({\Bbb D}, B_Y)$, $g(0) =x$ and $g(r) =y$.
It follows that
\[
K_Y (x,y) = K_Y (g(0), g(r)) \leq \rho(0,r) = K_{\mbox{lim}} + \varepsilon.
\]
Since this is true for all $\varepsilon > 0$, we have
\begin{equation}
\label{one-way}
K_Y (x,y) \leq K_{\mbox{lim}} = \lim_{\cal U} K_X ( x_i, y_i).
\end{equation}

To establish the reverse inequality, we use the fact that $K_X = C_X$
to select functions $f_i \in H(B_X, {\Bbb D})$ satisfying
$f_i(x_i) =0$, $f_i(y_i) > 0$ and $K_X(x_i,y_i) = \rho(0, f_i(y_i))$.
Applying part (i), we get a function $F \in H(B_Y, {\Bbb D})$
which satisfies $F(x) =0$ and
\[
\rho(0, F(y) ) = \lim_{\cal U} \rho(0, f_i(y_i)) = \lim_{\cal U}
K_X(x_i,y_i) = K_{\mbox{lim}}.
\]
But now the distance decreasing property of Kobayashi distances
under $F$ allows us to conclude that
\[
K_{\mbox{lim}} = \rho(0, F(y)) = \rho( F(x), F(y) ) \leq K_Y(x,y).
\]
Combining this with (\ref{one-way}) completes the proof.
\qed

Our next result is a new proof of a result of Davie and Gamelin~\cite{DG}.
It has come to our attention that M. Lindstr\"{o}m and R. Ryan have
independently obtained a proof of this result using ultrapower
techniques.

\begin{theorem}
\label{extend-f}
For $X$ a complex Banach space, there exists an algebra homomorphism
of norm one,
\[
E : H^\infty(B_X) \to H^\infty(B_{ X^{**} } )
\]
which satisfies
\[
Ef \big | {B_X} = f
\]
for all $ f \in H^\infty (B_X)$.
\end{theorem}

\noindent
{\sc Proof}: Given $X$ and
$f \in H^\infty (B_X)$, we choose an ultrafilter $\cal U$ according to the
principle of local reflexivity. To apply Theorem~\ref{extend-ultra}(i),
we take $f_i = f/ \| f \|_\infty$ for all $i$. Unless $f$ is constant,
the hypothesis (\ref{hypoth-limit}) is satisfied for $a_i =0$. In any
case, it follows that the function
\begin{eqnarray*}
F \colon B_Y &\to& \IC\\
F( (x_i)_{i \in I} + N_{\cal U} ) &=& \lim_{\cal U} f (x_i)
\end{eqnarray*}
is holomorphic and $\| F \|_\infty \leq \| f \|_\infty$.
Put $Ef= F \circ J$, where $J$ is as in the principle of local
reflexivity. It is straightforward to check that $E$ has the
required linearity and multiplicative properties and that
$\|E \| \leq 1$. Since $J$ coincides with the canonical embedding
$j \colon X \to (X)_{\cal U}$ on $B_X$, it is also easy to see that
$Ef$ coincides with $f$ on $B_X$. \qed

\begin{theorem}
\label{extend-k}
For $X$ a complex Banach space,
\[
K_{ X^{**} } (x,y) = K_X (x,y)
\]
for all $ x,y \in B_X$.
\end{theorem}

\noindent
{\sc Proof}:
Since the canonical inclusion from $B_X$ to $B_{X^{**}}$ is
continuous and linear, it is holomorphic and the distance decreasing
property of the Kobayashi metric implies
\[
K_{X^{**}} (x,y) \leq K_X (x,y)
\]
for $x,y \in B_X$.

Fix $x,y \in B_X$.
By Montel's theorem, we can find $f \in H(B_X, {\Bbb D})$ so that
$f(0) =0$ and
\[
K_X(x,y) = C_X(x,y) = \rho(f(x), f(y)).
\]
By Propostion~\ref{extend-f}, we can find an extension $\tilde{f} \in
H^\infty(B_{X^{**}})$ of $f$ with $\|\tilde{f}\|_\infty \leq 1$.
Since $\tilde{f}(x) = f(x) \in \Bbb D$, $\tilde{f}$ has all its
values in $\Bbb D$. We conclude
\[
K_X(x,y) = \rho( f(x), f(y)) = \rho( \tilde{f}(x), \tilde{f}(y) )
\leq K_{X^{**}} (x,y). \qed
\]

We are now in a position to extend Th\'{e}or\`{e}me~4.3 of \cite{DTV}.

\begin{theorem}
\label{1-compl}
If a complex Banach space $X$ is 1-complemented in its second dual,
then $B_X$ is C-connected.
\end{theorem}

\noindent
{\sc Proof}: Let $P$ denote a norm 1 projection from $X$** onto $X$
and let $p \neq q \in B_X$. Since $B_{X^{**}}$ is the unit ball of a dual
Banach space, Th\'{e}or\`{e}me~4.3 of \cite{DTV} implies that there exists
$\phi \in H( {\Bbb D}, B_{X^{**}} )$ and $u,v \in \Bbb D$
satisfying $\phi(u) =p$, $\phi(v) =q$ and $\rho(u,v) = K_{X^{**}} (p,q)$.

Now $P \circ \phi \in H( {\Bbb D}, B_X )$ and $P \circ \phi (u) = p$,
$P \circ \phi (v) = q$. By Theorem~\ref{extend-k} and
\cite[Th\'{e}or\`{e}me~2.5]{DTV},
\[
\rho (u,v) = K_{X^{**}} (p,q) = K_X (p,q) = C_X (p,q).
\]
Hence $P \circ \phi$ is a complex geodesic in $B_X$. Since $p$ and $q$
are arbitrary, $B_X$ is C-connected. \qed

\begin{remark}
Preduals of $C$*-algebras satisfy the hypotheses of Theorem~\ref{1-compl}
(see \cite{takesaki})
and these include examples which are not covered by the results in \cite{DTV}.

On the other hand, $c_0$ is well known not to satisfy the hypotheses of
Theorem~\ref{1-compl}, although $B_{c_0}$ is C-connected (which
follows from homogeneity --- see Remarks~\ref{l-infty}).
\end{remark}

\section{UNIQUENESS RESULTS}
\setcounter{equation}{0}
\label{sec-uniq}
If $\phi$ is a complex geodesic joining the points $p$ and $q$ of
a domain $\cal D$ and $f$ is a biholomorphic automorphism of
$\Bbb D$, then $\phi \circ f$ is also a complex geodesic joining
$p$ and $q$ (because $f$ is a $\rho$-isometry). Thus there is never
a unique complex geodesic joining $p$ and $q$, because
of this possibility of reparametrizing complex geodesics. However,
Vesentini~\cite{vesentini} has shown that if $\phi$ and $\psi$
are complex geodesics then they have the same range $\phi({\Bbb D})
= \psi( {\Bbb D})$ if and only if $\psi = \phi \circ f$ for some
biholomorphic automorphism $f$ of $\Bbb D$ (this can also be deduced
from the global vector-valued subordination theorem of
Finkelstein and Whitley~\cite{Fink-W}). We thus discuss uniqueness
of complex geodesics up to reparametrization, by
means of the following normalization.

We call a complex geodesic $\phi$ a {\em normalized
geodesic joining $p$ and $q$} if $\phi(0) = p$ and
$\phi(s) =q$ for some positive real number $s$. The number $s$
is uniquely determined by $p$ and $q$ --- in fact $s = \tanh C_{\cal D}(p,q)$.
By the homogeneity of the unit disc and the result of Vesentini
cited above, it follows that there is a unique normalized complex
geodesic joining two points $p, q \in {\cal D}$ if
and only if all complex geodesics joining $p$ and $q$
have the same range.

The following are known results concerning uniqueness.

\begin{description}
\item{(a)} If $B_X$ is the unit ball of a Banach space $X$ and
$x \in B_X$, $x \neq 0$, then there is a unique normalized complex
geodesic joining $0$ and $x$ if and only if $x/\|x\|$
is a complex extreme point of $B_X$. (Vesentini~\cite{vesentini}).

\item{(b)} If $\cal D$ is a strictly convex domain (i.e.~each
point of the boundary $\partial \cal D$ is a (real) extreme point
of $\cal D$) in a finite dimensional space, then there exist
unique normalized complex geodesics joining all
pairs of points in $\cal D$. (Lempert~\cite{lempert-bull}
).
\end{description}

In this section, we extend (b) to a class of Banach spaces
which includes all reflexive Banach spaces and give a general
criterion for uniqueness of complex geodesics
which highlights the problem of interpolating between
the results (a) and (b) above.
A more detailed study of non-uniqueness of complex geodesics has been
undertaken by
Gentili~\cite{gent-nu,gent-bal,gent-l1} (see also Section~\ref{sec-taut}).

\begin{definition}
\label{arnp-def}
A complex Banach space $X$ has the {\em analytic Radon-Nikodym property}
(aRNP) if each $f \in H^\infty ({\Bbb D}, X)$ has radial limits
almost everywhere on the unit circle.
\end{definition}

$H^\infty( {\Bbb D}, X)$ means those functions in $H ( {\Bbb D}, X)$
which have bounded range. If $X$ has aRNP and $f \in H^\infty ( {\Bbb D}, X)$,
we can extend $f$ to almost all points $e^{i\theta} \in \partial {\Bbb D}$
(almost all with respect to Lebesgue measure on $\partial \Bbb D$)
by
\[
f(e^{i \theta}) = \lim_{r \to 1^-} f (r e^{i\theta}).
\]
Moreover $f$ is uniquely determined by the boundary values
$f (e^{i\theta})$. Reflexive Banach spaces and Banach spaces with
the Radon-Nikodym property (e.g separable dual spaces) have aRNP.
The Banach space $c_0$ does not have aRNP.
For further details we refer to \cite[Chapter 12]{dineen} and
\cite{haagerup}.

Now suppose we have a convex bounded domain $\cal D$ in a Banach
space $X$ with aRNP. Let $p \neq q \in \cal D$ and let $G(p,q)$
denote the set of all normalized complex geodesics
joining $p$ and $q$. If $\phi, \psi \in G(p,q)$ and $0 < \lambda < 1$,
and $ s = \tanh C_{\cal D} (p,q)$, then $\lambda \phi + (1-\lambda)\psi
\in G(p,q)$. This follows from convexity of $\cal D$ (which ensures
that $\lambda \phi + (1-\lambda)\psi ({\Bbb D}) \subset {\cal D}$)
together with the facts that
\[
(\lambda \phi + (1-\lambda)\psi)(0) = \phi(0) = \psi(0) = p \mbox{ and }
(\lambda \phi + (1-\lambda)\psi)(s) = \phi(s) = \psi(s) = q.
\]
In other words, $G(p,q)$ is a convex subset of $H({\Bbb D}, {\cal D})$

Next, notice that if $\phi \in G(p,q)$, then
\[
\lim_{r \to 1^-} C_{\cal D} ( \phi (r e^{i\theta}), \phi(0))
= \lim_{r \to 1^-} \rho ( r e^{i\theta},0) = \infty.
\]
It follows that at points $e^{i\theta}$ where $\phi(e^{i\theta})$
is defined (almost all points on $\partial \Bbb D$ by aRNP),
$\phi(e^{i\theta}) \in \partial \cal D$. Now if $\psi$ is another
element of $G(p,q)$, then for almost all $\theta$,
\[
\phi(e^{i\theta}), \psi(e^{i\theta}) \mbox{ and }
\left( \frac{\phi+ \psi}{2}\right) (e^{i\theta}) =
(1/2) (\phi(e^{i\theta}) + \psi(e^{i\theta}))
\]
are all in $\partial \cal D$. If we now assume that
$\cal D$ is strictly convex, then we must have
$\phi(e^{i\theta}) = \psi(e^{i\theta})$ for almost all
$e^{i\theta}$. This implies $\phi=\psi$. We have thus
proved the following result.

\begin{theorem}
\label{uniq-arnp-th}
If $X$ is a complex Banach space with the analytic Radon-Nikodym
property and ${\cal D} \subset X$ is a strictly convex bounded
domain, then there exists at most one normalized complex
geodesic joining $p$ and $q$.
\end{theorem}

We now restrict our attention to the case where $\cal D$ is
the open unit ball $B_X$ of $X$. Let $\phi, \psi \in G(p,q)$
for two points $p,q \in B_X$ and suppose again that $X$ has
aRNP. Let $g = \psi - \phi$ and $s = \tanh C_X(p,q)$. Since
$\phi(0) = \psi(0) = p$ and $\phi(s) = \psi(s) =q$, $g(0) = g(s) = 0$
and we can therefore write $g(z) = z(z-s)h(z)$ for some
$h \in H^\infty({\Bbb D}, X)$. Using aRNP and the convexity
of $G(p,q)$, we can see that the following result holds.

\begin{proposition}
\label{uniq-equiv}
If $X$ has aRNP and $\phi$ is a normalized complex geodesic
joining two points $p,q \in B_X$, then $\phi$ is the unique
such geodesic if and only if the zero function is the only
element $h \in H^\infty({\Bbb D} ,X)$ satisfying
\begin{equation}
\label{h-condition}
\| \phi(e^{i\theta}) + \lambda e^{i\theta} (e^{i\theta}-s)h(e^{i\theta})\|
=1
\end{equation}
for almost all $\theta$, all $\lambda \in [0,1]$ (where $s =
\tanh C_X(p,q)$).
\end{proposition}

\begin{examples}
\begin{description}
\item{(a)} If $X$ has aRNP and is 1-complemented in $X^{**}$,
and if $B_X$ is strictly convex, then Theorems~\ref{1-compl}
and \ref{uniq-arnp-th} show that there exists a unique normalized
complex geodesic joining each pair of points $p,q \in B_X$.

\item{(b)} This applies in particular when $X$ is the space $\ell^p$
of $p$-summable sequences ($1 < p < \infty$), because $\ell^p$ is reflexive
and has a strictly convex unit ball.
\end{description}
\end{examples}
\section{CONTINUOUS COMPLEX GEODESICS}
\label{sec-cont}
\setcounter{equation}{0}

In this section we show that complex
geodesics can be extended continuously to the boundary
under a complex uniform convexity hypothesis.

\begin{definition}
\label{cuc-def}
If ${\cal D} \subset X$ is a domain in a complex Banach space $X$, then we
define
\[
\delta_{\cal D} (z,v) = \sup \{ r > 0 : z + rv {\blb D} \subset
{\cal D} \}
\]
for $z \in \cal D$, $v \in X$, $\|v\| = 1$.

We define the {\em modulus of complex convexity} of $\cal D$ to be
\[
\delta_{\cal D} (\varepsilon) = \sup \{ \delta_{\cal D} ( z, v) : z \in {\cal
D}, d(z, \partial {\cal D}) \leq \varepsilon, \|v\| = 1 \},
\]
(where $d(z, \partial {\cal D} )$ denotes the distance from a point $z$ in
$\cal D$ to the boundary).

The domain $\cal D$ is called {\em complex uniformly convex} if
$\delta_{\cal D}( \varepsilon) \to 0$ as $\varepsilon \to 0$.
\end{definition}

\begin{remarks}
For the case where ${\cal D} = B_X$ is the unit ball of
$X$, Globevnik~\cite{globevnik} introduced
\[
\omega_c(\varepsilon) = \sup \{ \| y\| : \|x + \zeta y\| \leq 1 +
\varepsilon \mbox{ for all } \|x\| = 1, \zeta \in {\blb D} \}
\]
and this function is closely related to $\delta_{\cal D} ( \varepsilon )$.
In fact, it is easy to check that
\[
\frac{\omega_c(\varepsilon)}{1 + \varepsilon} =
\sup\{ \delta_{\cal D} (z,v) : \|v\| = 1, \|z\| = 1/(1+\varepsilon) \} =
\delta_{\cal D} \left( \frac{\varepsilon}{1+\varepsilon} \right)
\]
and thus that
\[
\delta_{\cal D}( \varepsilon/2) \leq \omega_c(\varepsilon) \leq 2
\delta_{\cal D} (\varepsilon) \quad (0 < \varepsilon \leq 1).
\]

Functions which are inverse to $\delta_{\cal D}(\varepsilon)$ and
$\omega_c( \varepsilon)$ were considered by Davis, Garling and
Tomczak-Jaegermann~\cite{DGT} and called $h^X_\infty$ and $H^X_\infty$
(respectively). Dilworth~\cite[theorem 2.1]{dilworth}
has shown that complex uniform
convexity of $B_X$ (or uniform $H_\infty$-convexity of $X$ in the
notation of \cite{DGT}) is equivalent to the notion of uniform
$PL$-convexity which was studied intensively in \cite{DGT}.
\end{remarks}

A result similar to the following one can be obtained for the case where
the domain is the unit ball using theorem~2 of \cite{globevnik}.
\begin{proposition}
\label{carath-cuc}
If ${\cal D} \subset X$ is a convex domain, then
\[
c_{\cal D}(z,v) \geq \frac{\|v\|}{2 \delta_{\cal D} (\varepsilon)}
\]
holds for $v \in X$, $z \in \cal D$, $\varepsilon = d(z, \partial{\cal
D})$.
\end{proposition}

\noindent
{\sc Proof}: We will make use of the fact that $c_{\cal D} = k_{\cal D}$
for $\cal D$ convex (see \cite{DTV}). Fix $z \in \cal D$, $v \in X$ and
consider holomorphic mappings $f \colon {\Bbb D} \to \cal D$ with $f(0) =
z$, $f'(0)=v/r$, $r > 0$.

Consider the function 
\[
g(\zeta) = \frac{1}{2\pi} \int_0^{2\pi} f( e^{i\theta} \zeta)
(1 + \cos \theta) \, d\theta.
\]
Since $(1+  \cos \theta)\, d\theta/(2\pi)$ is a probability measure on
$[0,2\pi]$, we may regard $g(\zeta)$ as a limit of convex combinations
of points in $\cal D$. So $g$ takes values in the closure $\bar{\cal D}$.
Since $g(0) = f(0) \in \cal D$ and $\cal D$ is convex, it follows
that $g$ maps $\Bbb D$ into $\cal D$ \cite[p.~376]{vesentini}.

Working with the power series representation of $f$, we find that
\[
g(\zeta) = f(0) + \frac{1}{2} f'(0) \zeta.
\]
(In fact, $g(\zeta)$ is a Cesaro mean of the power series
of $f$.). Thus
\[
g(\zeta) = f(0) + \frac{v}{2r} \zeta
= z + \frac{\|v\|}{2r} \frac{v}{\|v\|} \zeta
\]
maps $\Bbb D$ into $\cal D$, which shows that
\[
\frac{\|v\|}{2r} \leq \delta_{\cal D}( z , v/\|v\| ) \leq \delta_{\cal D}
(\varepsilon).
\]
Rearranging this, we find $ r \geq \|v\| / (2 \delta_{\cal
D}(\varepsilon))$. Since $k_{\cal D}(z,v) = c_{\cal D}(z,v)$
is the infimum of all possible
values of $r$, the result follows. \qed

\begin{theorem}
\label{th-cont}
Let $\cal D$ be the unit ball $B_X$ of a complex Banach space $X$. If $\cal
D$  is complex uniformly convex
and $\delta_{\cal D}( \varepsilon) \leq A\varepsilon^s$ for some constants
$A >0$, $s>0$, then all complex geodesics $\phi\colon {\Bbb D} \to {\cal
D}$ extend to continuous functions
$\phi \colon \bar{\Bbb D} \to \bar{\cal D}$.
\end{theorem}

\noindent
{\sc Proof}: By Proposition \ref{carath-cuc}, we have (for $\zeta \in \Bbb D$)
\[
c_{\cal D} ( \phi(\zeta), \phi'(\zeta) ) \geq \frac{\|\phi'(\zeta)\|}
{2 \delta_{\cal D} ( 1 - \|\phi(\zeta)\| ) }.
\]
Using the hypothesis and the fact that $\phi$ is a complex geodesic,
we deduce
\[
c_{\cal D} (\phi(\zeta),\phi'(\zeta)) =
c_{\Bbb D} (\zeta,1) = \frac{1}{1 - |\zeta|^2} \geq \frac{\|\phi'(\zeta)\|}
{2A(1 - \|\phi(\zeta)\|)^s}
\]
or
\begin{equation}
\label{geod-est}
\|\phi'(\zeta)\| \leq \frac{2A(1 - \|\phi(\zeta)\|)^s}{1 - |\zeta|^2}.
\end{equation}

Next observe that
\[
C_{\cal D} ( \phi(0), \phi(\zeta) ) = \rho(0,\zeta) = \tanh^{-1}(|\zeta|)
\leq C_{\cal D} (0, \phi(\zeta) ) + C_{\cal D} (0, \phi(0)).
\]
Recall that $C_{\cal D} (0,z) = \tanh^{-1} \|z\|$ for $z \in {\cal D} = B_X$.
Using elementary estimates, we
conclude that
\[
1 - \|\phi(\zeta)\| \leq A_\phi (1-|\zeta|)
\]
where $A_\phi$ is a constant depending on $\|\phi(0)\|$.

Combining this observation with (\ref{geod-est}), we see that
\[
\|\phi'(\zeta)\| \leq A_1 \frac{(1-|\zeta|)^s}{ 1 - |\zeta|^2}
\leq A_1 \frac{1}{(1-|\zeta|)^{1-s} }
\]
which implies (see for instance \cite[theorem~5.5]{duren})
that $\phi$ satisfies a Lipschitz condition
\[
\| \phi(\zeta) - \phi(\eta) \| \leq C |\zeta - \eta|^s.
\]
Hence $\phi \colon {\Bbb D} \to X$ is uniformly continuous and
extends continuously to a function $f \colon \bar{\Bbb D} \to X$.
\qed

If a complex geodesic $\phi \colon {\Bbb D} \to {\cal D}$ extends to a
continuous function $\phi \colon \bar{\Bbb D} \to \bar{\cal D}$, we call
$\phi$ a {\em continuous complex geodesic}.

\begin{examples} 
\label{uniq-ex}
\begin{description}
\item{(i)} For $X = L^1$, ${\cal D} = B_X$, Globevnik~\cite{globevnik}
proved that $\delta_{\cal D} \leq A \sqrt{\varepsilon}$ (or rather, he
proves the equivalent fact that $\omega_c(\varepsilon) \leq A
\sqrt{\varepsilon}$). Thus all complex geodesics in $B_X$ are continuous,
by Theorem~\ref{th-cont}.

\item{(ii)} More generally, if $X$ is the predual of a $C$*-algebra and
${\cal D} = B_X$, we deduce from a result due to Haagerup (see
\cite[theorem 4.3]{DGT}) that $\delta_{\cal D}(\varepsilon) \leq A
\sqrt{\varepsilon}$.

Thus all complex geodesics in $B_X$ are continuous
by Theorem~\ref{th-cont}. Existence of complex geodesics in $B_X$
is guaranteed by Theorem~\ref{1-compl}.
\end{description}
\end{examples}

\begin{remark}
From corollary 2.5 of \cite{gent-bal}, it follows that if all complex
geodesics in the unit ball $B_X$ are continuous, then all points of
$\partial B_X$ are complex extreme points. Theorem~\ref{th-cont} falls short
of being a converse to this.
\end{remark}

\def\i{i}
\section{EXAMPLES IN CLASSICAL BANACH SPACES}
\setcounter{equation}{0}
\label{sec-ex}

We apply the results of the preceding sections to give a complete
description of the complex geodesics in the unit ball of $\ell^p$, $1 \leq
p < \infty$.
To obtain these examples, we require a Banach space version of a result of
Lempert \cite[proposition 1]{lempert-bull} and \cite{Roy-W}. Various
extensions are possible and we have chosen one which is suitable for the
applications we have in mind.

\begin{definition}
\label{hinfty-def}
For $X$ a complex Banach space with dual space $X^*$, we let
$H^\infty_*({\Bbb D}, X^*)$ denote the space of $X^*$-valued bounded analytic
functions on $\Bbb D$ which have weak*-radial limits at almost all boundary
points.

In other words $f \in H^\infty_*({\Bbb D}, X^*)$ means that $f \in
H^\infty({\Bbb D}, X^*)$
and there exists a function
$\tilde{f} \colon \partial {\Bbb D} \to X^*$ so that
\[
\lim_{r \to 1^-} \langle x, f(r e^{i\theta} )  \rangle =
\langle x, \tilde{f} ( e^{i\theta} ) \rangle \qquad (\mbox{all } x \in X)
\]
holds for almost all $\theta \in \IR$.
\end{definition}

If $X^*$ has aRNP, then all functions in $H^\infty({\Bbb D}, X^*)$
have norm radial
limits at almost all points of $\partial \Bbb D$, and therefore
$H^\infty_* ({\Bbb D},  X^*) = H^\infty({\Bbb D}, X^*)$.
By a result of Danilevich~\cite[theorem 1.4]{danilevich}, this
equality also holds if $X$ is
separable. 

In general, the limit function
$\tilde{f}(e^{i\theta})$ may only be weak*-measurable, a rather
intractable condition. Moreover it is possible that the space
$H^\infty_*({\Bbb D}, X^*)$ depends on the choice of a
predual $X$ for $X^*$.
However, the function $\tilde{f}( e^{i\theta})$
determines the holomorphic function $f(\zeta)$ uniquely as can be seen by
applying
standard results (see \cite{duren}) to the scalar-valued bounded analytic
functions $\langle x, f(\zeta) \rangle$. We will therefore not cause
confusion by using the notation $f(e^{i\theta})$ instead of
$\tilde{f} (e^{i\theta})$ for the boundary function.

\begin{lemma}
\label{radial}
Let $X$ be a complex Banach space, $f \in H^\infty_*({\Bbb D}, X^*)$ and
$h \colon \bar{\Bbb D} \to X$ a continuous function which is holomorphic on
$\Bbb D$. Then
\[
\lim_{r \to r^-} \langle h(r e^{i\theta}), f(r e^{i\theta}) \rangle
= \langle h(e^{i\theta}), f(e^{i\theta}) \rangle
\]
for almost all $\theta$ and
\[
\langle h(\zeta),  f(\zeta) \rangle
= \frac{1}{2\pi} \int_0^{2\pi}
\langle h(e^{i\theta}),  f(e^{i\theta}) \rangle
\frac{1 - |\zeta|^2}{1 + |\zeta|^2 - 2 \mbox{\rm Re}\, (e^{-i\theta}\zeta) } \,
d\theta
\]
for all $\zeta \in \Bbb D$.
\end{lemma}

\noindent
{\sc Proof}: The first assertion follows from the inequality
\begin{eqnarray*}
\lefteqn{|\langle h(r e^{i\theta}), f(r e^{i\theta}) \rangle
- \langle h(e^{i\theta}), f(e^{i\theta}) \rangle|}\\
&\leq& |\langle h(r e^{i\theta}) - h(e^{i\theta}), f(r e^{i\theta}) \rangle|
+ |\langle h(e^{i\theta}), f(r e^{i\theta}) - f(e^{i\theta}) \rangle|\\
&\leq& \|h(r e^{i\theta}) - h(e^{i\theta})\| \, \|f\|_\infty 
+ |\langle h(e^{i\theta}), f(r e^{i\theta}) - f(e^{i\theta}) \rangle|\\
& \to & 0
\end{eqnarray*}
for almost all $\theta$ by continuity of $h$ and the definition of
$H^\infty_*({\Bbb D}, X^*)$. The second assertion follows from the standard
fact (see \cite{duren}) that scalar-valued bounded analytic functions
like $\langle f(\zeta), h(\zeta) \rangle$ are the Poisson integrals
of their (almost everywhere) boundary values. \qed

If $X$ is a Banach space and $x \in \partial B_X$, then the Hahn-Banach
theorem assures us of the existence of at least one supporting hyperplane
for $B_X$ at $x$. That is, there exists $N_x \in X^*$ such that
$\langle x, N_x \rangle = 1$ and $\mbox{Re}\,\langle p,  N_x \rangle < 1$
for $p \in B_X$. We will use the
notation $N_x$ for a choice of one such functional, bearing in mind the
possibility that it may not be unique.

\begin{proposition}
\label{lempert-criterion}
Let $X$ be a complex Banach space and
$\phi \colon \bar{\Bbb D} \to \bar {B}_X$ a continuous map satisfying
\begin{description}
\item{(i)} $\phi \mid {\Bbb D}$ is holomorphic and $\phi({\Bbb D}) \subset
B_X$;
\item{(ii)} $\phi( \partial {\Bbb D}) \subset \partial B_X$; and
\item{(iii)} there exists a choice of $N_{\phi(\zeta)}$ for almost all $\zeta \in
\partial \Bbb D$ and a measurable function $p \colon \partial{\Bbb D} \to
\IR^+$ such that the mapping
\[
h(e^{i\theta}) =
e^{i\theta} p(e^{i\theta}) N_{\phi(e^{i\theta})}
\]
is almost everywhere the weak*-radial limit of a function
$h \in H^\infty_*({\Bbb D}, X^*)$.
\end{description}
Then $\phi$ is a complex geodesic.
\end{proposition}

\noindent
{\sc Proof}: Let $g \colon {\Bbb D} \to B_X$ be a holomorphic mapping with
$g(0) = \phi(0)$ and $g'(0) = \lambda \phi'(0)$, $\lambda \geq 0$. Let
$g_r(\zeta) = g (r \zeta)$ for $0 < r < 1$ and $\zeta \in {\Bbb D}$. Then
$g_r$ is continuous on $\bar{\Bbb D}$ and holomorphic on $\Bbb D$, $g_r(0)
= \phi(0)$ and $g_r'(0) = \lambda r \phi'(0)$. Moreover $g_r(\zeta) \in
B_X$ for $\zeta \in \bar{\Bbb D}$. From the hypotheses, we see that
\[
1 = \langle \phi(e^{i\theta})  ,  N_{\phi(e^{i\theta})}   \rangle >
\mbox{Re}\, \langle g_r(e^{i\theta}), N_{\phi(e^{i\theta})} \rangle
\]
for almost all $\theta$.
Hence
\[
\mbox{Re\, }\left\langle
 \frac{\phi(e^{i\theta}) - g_r(e^{i\theta})}{e^{i\theta}} ,
e^{i\theta} p(e^{i\theta}) N_{\phi(e^{i\theta})}
\right\rangle
= p(e^{i\theta})  \mbox{Re}\,
\langle \phi(e^{i\theta}) - g_r(e^{i\theta}), N_{\phi(e^{i\theta})} \rangle
> 0
\]
for almost all $ \theta$.

Since $\displaystyle\frac{\phi(\zeta) - g_r(\zeta)}{\zeta}$ is
holomorphic on $\Bbb D$
and continuous on $\bar{\Bbb D}$, the function
\[
H(\zeta) = \left\langle \frac{\phi(\zeta) - g_r(\zeta)}{\zeta}, h(\zeta)
\right\rangle
\]
is the Poisson integral of its boundary values $H(e^{i\theta})$
by Lemma~\ref{radial}. By the above remark and Lemma~\ref{radial},
$\mbox{Re}\, H(e^{i\theta}) > 0$ for almost all $\theta$ and it follows from
the Poisson formula that $\mbox{Re}\,H(0) > 0$, {\em i.e.}
\[
\mbox{Re}\, \langle \phi'(0) - g_r'(0), h(0) \rangle =
(1-\lambda r) \mbox{Re}\, \langle \phi'(0) , h(0) \rangle > 0.
\]
Applying this to the special case where $g(\zeta) = \phi(0)$ is constant
(and $\lambda =0$) we see that
$\mbox{Re}\, \langle \phi'(0), h(0) \rangle  > 0$.
Thus, returning to the general case, we have $1-r\lambda > 0$. Since this
is true for all $0 < r < 1$, we deduce that $\lambda \leq 1$.

Since this is true for all $g$, we have established
\[
k_X(\phi(0), \phi'(0)) = c_X(\phi(0), \phi'(0)) = 1
\]
which shows, by Proposition~\ref{geod-equiv},
that $\phi$ is a complex geodesic. \qed

For $\mu$ a $\sigma$-finite measure on a $\sigma$-algebra $\Sigma$ of
subsets of a set $\Omega$ and $1 \leq p < \infty$, we use the standard
notation $L^p(\mu)$ for the Banach space of (equivalence classes of)
$p$-summable $\Sigma$-measurable ($\IC$-valued)
functions on $\Omega$ normed by
$\|f\|_p = \left( \int_\Omega |f|^p \, d\mu \right)^{1/p}$. $L^\infty(\mu)$
denotes the essentially bounded $\Sigma$-measurable functions with
the essential sup norm $\|f\|_\infty$. These include as special cases
the sequence spaces $\ell^p$ (where $\mu$ is counting measure on the
natural numbers) and the finite-dimensional spaces $\ell^p_n$
(which are $\IC^n$ with the norm
$\|(z_1, z_2, \ldots, z_n)\|_p = \left(\sum_{j=1}^n |z_j|^p \right)^{1/p}$
for $1 \leq p < \infty$).

In order to discuss complex geodesics in the unit ball $B_p$ of
$L^p(\mu)$ for $1 \leq p < \infty$ we will consider nonconstant mappings
$\phi \colon {\Bbb D} \to B_p$ of the form
\begin{equation}
\label{geod-form}
\phi(\zeta)(\omega) = c(\omega)
\left( \frac{\zeta - \alpha(\omega)}{1 - \overline{\alpha(\omega)} \zeta}
\right)^{\beta(\omega)}
\left( \frac{1 - \overline{\alpha(\omega)} \zeta}{1 - \bar{\gamma} \zeta}
\right)^{2/p}
\end{equation}
($\zeta \in \Bbb D$, $\omega \in \Omega$)
where the parameter $\gamma$ and the measurable functions
$\alpha(\omega)$, $\beta(\omega)$ and $c(\omega)$ satisfy
\begin{equation}
\label{param-conds}
\left.
\begin{array}{ll}
(a) & \gamma \in {\Bbb D}, \quad \alpha \in L^\infty(\mu), \quad
                                  \|\alpha\|_\infty \leq 1, \\
    & \beta \mbox{ takes only the values } 0 \mbox{ and } 1.\\[1mm]
(b) & \int_\Omega |c(\omega)|^p (1 + |\alpha(\zeta)|^2) \, d\mu(\omega)
                                    = 1 + |\gamma|^2\\[1mm]
(c) & \int_\Omega |c(\omega)|^p \alpha(\omega) \, d\mu = \gamma
\end{array}
\right\}
\end{equation}
Later, we will specialise to the case where $L^p(\mu) = \ell^p$ and then
we will start to use subscript notation --- $\phi_j(\zeta)$ rather than
$\phi(\zeta)(j)$, $\alpha_j$ instead of $\alpha(j)$, {\em etc.} --- and of
course summation over $j$ in place of integrals.

\begin{proposition}
\label{Lp-geod}
Let $B_p$ denote the open unit ball of $L^p(\mu)$, $1 \leq p < \infty$.
Then every nonconstant mapping $\phi$ of the form (\ref{geod-form})
where $\gamma$, $\alpha(\omega)$, $\beta(\omega)$ and $c(\omega)$
satisfy the conditions (\ref{param-conds}) is a complex geodesic in $B_p$.
\end{proposition}

Note that at points where $c(\omega) = 0$, the values of $\alpha(\omega)$
and $\beta(\omega)$ are immaterial. Thus we can suppose if we
wish that
\[
\mbox{support}\,(\beta) = \{ \omega: \beta(\omega) = 1 \} \subset
\{ \omega : c(\omega) \neq 0 \} = \mbox{support}\,(c)
\]
and that $\mbox{support}\,(\alpha) \subset \mbox{support}\,(c)$. Since
$(\zeta - \alpha(\omega))/(1 - \overline{\alpha(\omega)}\zeta)$ is a
constant function of $\zeta$ when $|\alpha(\omega)| =1$ (the constant is of
modulus $1$), we can also
assume $\mbox{support}\,(\beta) \subset \{ \omega : |\alpha(\omega)| < 1 \}$.
(Then there is no problem defining $\phi(\zeta)$ for $|\zeta| = 1$.)
With these assumptions, the condition for $\phi$ to be non-constant
is
\[
\mu( \{ \omega : c(\omega) \neq 0 \} \cap \{ \omega : \beta(\omega) = 1
\mbox{ or } \alpha(\omega) \neq \gamma \} ) > 0.
\]
Note also that condition (\ref{param-conds})(b) implies that $c \in L^p(\mu)$.

\noindent
{\sc Proof}: Observe that for
$\delta = \left( \displaystyle\frac{2}{1 - |\gamma|} \right)^{2/p}$,
\begin{equation}
\label{delta-estimate}
|\phi(\zeta)(\omega)| \leq \delta |c(\omega)|
\qquad (\zeta \in \bar{\Bbb D}, \omega \in \Omega)
\end{equation}
so that $\phi(\bar{\Bbb D}) \subset L^p(\mu)$. Next if $(\zeta_n)$
is a sequence in $\bar{\Bbb D}$ which converges to a point
$\zeta \in \bar{\Bbb D}$, then $\phi(\zeta_n)(\omega) \to \phi(\zeta)(\omega)$
for all $\omega$. Hence by the Lebesgue dominated convergence theorem and
(\ref{delta-estimate}), $\| \phi(\zeta_n) - \phi(\zeta)\|_p \to 0$
as $n \to \infty$. This shows that $\phi \colon \bar{\Bbb D} \to L^p(\mu)$
is a continuous mapping.

For each fixed $\omega \in \Omega$ the map $\zeta \mapsto
\phi(\zeta)(\omega)$ is analytic on $\Bbb D$ and from (\ref{delta-estimate})
we see that its Taylor series expansion
\begin{equation}
\label{taylor}
\phi(\zeta)(\omega) = \sum_n a_n(\omega) \zeta^n
\end{equation}
with coefficients satisfying $|a_n(\omega)| \leq \delta |c(\omega)|$
(all $n$ and all $\omega$). Calculating the $a_n(\omega)$ using the
binomial theorem and multiplication of power series, we can
check that the functions $a_n(\omega)$ are measurable and then the estimate
on the coefficients implies that $a_n \in L^p(\mu)$.

Now, for a fixed $\zeta \in \Bbb D$ the sequence
$\sum_{j=0}^n a_j( \cdot ) \zeta^j$ of measurable functions converges
pointwise to $\phi(\zeta)(\cdot)$ as $n \to \infty$ by (\ref{taylor}).
Since
\[
\left| \sum_{j=0}^n a_j( \omega ) \zeta^j \right|
\leq \delta |c(\omega)| \sum_{j=0}^\infty |\zeta|^j =
\frac{\delta |c(\omega)|}{1 - |\zeta|}
\]
and $\delta c(\cdot)/(1-|\zeta|) \in L^p(\mu)$, the Lebesgue dominated
convergence theorem shows that
\[
\phi(\zeta)(\cdot) = \sum_{n=0}^\infty a_n(\cdot) \zeta^n
\]
in $L^p(\mu)$ for each $\zeta \in \Bbb D$. Hence $\phi$ is holomorphic on
$\Bbb D$.

For $\theta \in \IR$, we have
\begin{eqnarray*}
\|\phi(e^{i\theta})\|_p^p &=& \int_\Omega |\phi(e^{i\theta})(\omega)|^p
                                   \, d\mu(\omega)\\
&=& \int_\Omega |c(\omega)|^p \left|
           \frac{1- \overline{\alpha(\omega)}e^{i\theta} }
                {1 - \bar{\gamma}e^{i\theta} } \right|^2 \, d\mu(\omega)\\
&=& \frac{1}{1+|\gamma|^2 - 2 \mbox{Re}\,(\bar{\gamma}e^{i\theta}) }
         \int_\Omega |c(\omega)|^p (1 + |\alpha(\omega)|^2
                -2 \mbox{Re}\, (\overline{\alpha(\omega)}e^{i\theta})\, ) \,
                   d\mu(\omega)\\
&=& 1
\end{eqnarray*}
by conditions (\ref{param-conds}) (b) and (c). Hence $\phi(\partial
{\Bbb D}) \subset \partial B_p$ and $\phi( \bar{\Bbb D} ) \subset
\bar{B}_p$. Since $\phi$ is non-constant
and all unit vectors in $L^p(\mu)$ are complex extreme points of the
unit ball the strong maximum modulus theorem of Thorp and  Whitley
(see for instance \cite[proposition 6.19]{dineen}) implies that
$\phi({\Bbb D}) \subset B_p$.

We consider the dual space $(L^p(\mu))^*$ of $L^p(\mu)$ to
be identified with $L^q(\mu)$ ($1/p + 1/q = 1$, $1 < q \leq \infty$)
in a complex linear fashion (rather than the conjugate linear
identification frequently used) so that $g \in L^q(\mu)$ acts
on $f \in L^p(\mu)$ via
$\langle f,g \rangle = \int_\Omega f(\omega)g(\omega) \, d\mu(\omega)$.
If $f \in \partial B_p$, then one choice of a supporting hyperplane
$N_f \in L^q(\mu)$ is given by
\[
N_f(\omega) = {|f(\omega)|}^{p-2} \overline{f(\omega)}
\]
(where $0^{p-2} = 0$ for all $p$).

To complete the proof, we apply Proposition~\ref{lempert-criterion}
with
\begin{eqnarray*}
p (e^{i\theta}) &=& |1- \bar{\gamma} e^{i\theta}|^2\\
h(\zeta)(\omega)  &=& \tilde{c}(\omega)
                         (1 - \overline{\alpha(\omega)} \zeta)^{2-2/p}
                         (1 - \bar{\gamma} \zeta)^{2/p}
                         \left( \displaystyle\frac{\zeta - \alpha(\omega)}
                         {1 - \overline{\alpha(\omega)} \zeta}
                                               \right)^{1-\beta(\omega)}
\end{eqnarray*}
with $\tilde{c}(\omega) = |c(\omega)|^{p-2} \overline{ c(\omega) }$.
Note that $\tilde{c} \in L^q(\mu)$. If $1 < p < \infty$ an argument
similar to the one given above to show that $\phi$ is holomorphic
on $\Bbb D$ and continuous on $\bar{\Bbb D}$ shows that the same is
true of $h$. This is more than enough to show that
$h \in H^\infty_*({\Bbb D}, (L^p(\mu))^*)$. Since $h(e^{i\theta}) = e^{i\theta}
p(e^{i\theta}) N_{\phi(e^{i\theta})}$, this shows that $\phi$ is
a complex geodesic when $p > 1$.

When $p=1$,
\[
h(\zeta)(\omega) = (1 - \bar{\gamma} \zeta)^2 \tilde{c}(\omega)
                     \left( \displaystyle\frac{\zeta - \alpha(\omega)}
                         {1 - \overline{\alpha(\omega)} \zeta}
                                               \right)^{1-\beta(\omega)}
\]
and $\|\tilde{c}\|_\infty \leq 1$. It is quite easy to see that $h$
is holomorphic on $\Bbb D$ and that
$\|h(\zeta)\|_\infty \leq (1 + |\gamma|)^2$ for all $\zeta \in \Bbb D$.
For each $\omega \in \Omega$ and $\theta \in \IR$,
$\lim_{r \to 1^-} h(r e^{i\theta})(\omega) = h(e^{i\theta})(\omega)$.
Using the boundedness of $h$ and the Lebesgue dominated convergence theorem
it is then easy to see that $h$ has weak* radial limits at all points
$e^{i\theta} \in \partial \Bbb D$, {\em i.e.} that
\[
\lim_{r \to 1^-} \langle f, h(r e^{i\theta}) \rangle =
\langle f,  h(e^{i\theta}) \rangle
\]
for all $f \in L^1(\mu)$. Thus $h \in H^\infty_* ({\Bbb D},  (L^1(\mu))^* )$
and $\phi$ is also a complex geodesic in the $p =1$ case.
\qed

\begin{proposition}
\label{lp-exist}
Let $B_p$ denote the open unit ball of $\ell^p$, $1 \leq p < \infty$.
Then any two distinct points in $B_p$ can be joined by a
normalized complex geodesic $\phi$ of the form given in
(\ref{geod-form}).

That is, there exists a complex geodesic $\phi(\zeta) =
\big( \phi_j(\zeta) \big)_{j=1}^\infty$ joining the points where
\begin{equation}
\label{geod-form-N}
\phi_j(\zeta) =
c_j \left(\frac{\zeta - \alpha_j}
{ 1 - \bar{\alpha}_j \zeta}\right)^{\beta_j}
\left( \frac{ 1 - \bar{\alpha}_j \zeta}
{1- \bar{\gamma}\zeta} \right)^{2/p}
\end{equation}
$\gamma \in \Bbb D$, $\alpha_j \in \bar{\Bbb D}$, $c_j \in \IC$,
$\beta_j =0 \mbox{ or } 1$,
$ \sum_j |c_j|^p (1 +
|\alpha_j|^2) = 1+|\gamma|^2$ and
$ \sum_j |c_j|^p \alpha_j = \gamma$.
\end{proposition}

\noindent
{\sc Proof}:
The existence of some complex geodesic joining the two points
follows from Theorem~\ref{1-compl} (or
from \cite{DTV}). We also know that all complex geodesics in $B_p$
are continuous from Theorem~\ref{th-cont}. For the $p=1$ case the
uniform complex
convexity hypothesis was established by Globevnik \cite{globevnik} (see
Examples~\ref{uniq-ex}) and the $p > 1$ case is more straightforward because
$\ell^p$ is uniformly convex in the real sense (see \cite{DGT}).

Our proof that there exists geodesics of the
required form relies heavily on finite-dimensional results
({\em i.e.} results for the case of $\ell^p_n$)
of Poletski\v{\i}~\cite{poletskii} and
Gentili~\cite{gent-l1}.
For $1< p < \infty$, Poletski\v{i}~\cite{poletskii} proved that all
geodesics in the unit ball $B_{p,n}$ of $\ell^p_n$ are of the above form
(except that he omitted the possibility that $\beta_j$ could be $0$).
Gentili~\cite[Theorem 6]{gent-l1} proved that all continuous complex
geodesics in $B_{1,n}$ are of the above form. Now that
we know all complex geodesics in $B_{1,n}$ are continuous, it
follows that all complex geodesics in $B_{1,n}$
are of the above form. (We could actually circumvent Gentili's
result. By taking a limiting argument based on Poletski\v{\i}'s
result and the fact that $\cap_{p>1} B_{p,n} = B_{1,n}$, we could
establish the existence of complex geodesics of the required form
joining any pair of points in $B_{1,n}$.)

Now let $x,y \in B_p$ be two distinct points.

We consider $\ell^p_n$ as being identical with the subspace $\{(z_1, z_2,
\ldots, z_n, 0, 0, \ldots) \}$ of $\ell_p$ and we will use the notation
$x^{(n)}$ for the natural projection $(x_1, x_2, \ldots, x_n, 0, \ldots)$
of $x = (x_1, x_2, \ldots)$ into $\ell^p_n$. Similarly for $y^{(n)}$.
Consideration of the inclusion map $\colon B_{\ell^p_n} \to B_p$ and the
projection $z \mapsto z^{(n)}\colon B_p \to B_{\ell^p_n}$, which are both
holomorphic and therefore contractions with respect to the Carath\'{e}odory
distance, shows that
\[
C_{\ell^p_n} ( x^{(n)}, y^{(n)} ) = C_{\ell^p} ( x^{(n)}, y^{(n)} )
\leq C_{\ell^p} ( x, y ) \quad (x,y \in B_p).
\]

Let
$\phi^{(n)}$ denote a normalized complex geodesic in the unit ball of
$\ell^p_n$ with $\phi^{(n)}(0) = x^{(n)}$ and $\phi^{(n)}(s_n) = y^{(n)}$,
where $s_n = \tanh^{-1} C_{\ell^p}(x^{(n)}, y^{(n)}) > 0$. (Strictly
speaking this may not make sense for small $n$ when it may happen that
$x^{(n)} = y^{(n)}$.) Let the parameters associated with $\phi^{(n)}$ be
denoted $\alpha_j^{(n)}$, $\beta_j^{(n)}$, $\gamma^{(n)}$ and $c_j^{(n)}$
($1\leq j \leq n$).

There is a subsequence of $n$'s along which we have
\[
\alpha_j^{(n)} \to \alpha, \quad
\beta_j^{(n)} \to \beta_j, \quad
\gamma^{(n)} \to \gamma, \quad
c_j^{(n)} \to c_j
\]
for each $j$. Since $\sum_j |c_j^{(n)}|^p < 2$ for all $n$,
we have $\sum_j |c_j|^p
\leq 2$.
We claim that
$|\gamma| < 1$. Observe that
\[
\phi_j^{(n)}(\zeta) \to \phi_j(\zeta) =
c_j \left(\frac{\zeta - \alpha_j}
{ 1 - \bar{\alpha}_j \zeta}\right)^{\beta_j}
\left( \frac{ 1 - \bar{\alpha}_j \zeta}
{1- \bar{\gamma}\zeta} \right)^{2/p}
\]
uniformly on compact subsets of $\Bbb D$ as $ n \to \infty$ along the
subsequence. Hence $\phi_j(0) = x_j$ and $\phi_j(s) = y_j$ where
\[
s = \lim_n s_n = \lim_n \tanh^{-1} C_{\ell^p} (x^{(n)}, y^{(n)} )
= \tanh^{-1} C_{\ell^p} (x,y).
\]
Thus we can pick a $j$ with $\phi_j$ non-constant. If $|\gamma| = 1$,
then $\phi_j$ would be unbounded on $\Bbb D$ unless $\alpha_j = \gamma$.
$\phi_j$ unbounded leads to a
contradiction since each $\phi_j^{(n)}$ is bounded by $1$ and in the
case $\alpha_j = \gamma$,  $\phi_j$ would be constant. Hence $|\gamma| < 1$.
Using $\sum_j |c_j|^p < \infty$ we can argue as in the proof of
Proposition~\ref{Lp-geod} to show that
\[
\phi = (\phi_1, \phi_2, \ldots) \colon {\Bbb D} \to \ell^p
\]
is holomorphic on $\Bbb D$ and continuous on $\bar{\Bbb D}$.

Since $\sum_j |\phi_j^{(n)}(\zeta)|^p \leq 1$ for all $n$,
we have $\sum_j |\phi_j(\zeta)|^p \leq 1$ for
$|\zeta| < 1$.
Since $\phi(0) = x$ we have
$\phi({\Bbb D}) \subset B_p$ (see \cite[p.~376]{vesentini})
and since $\phi(s) = y$,
$\phi$ must be a complex geodesic.
It follows that $\|\phi(\zeta)\|_p^p = 1$ for $|\zeta| = 1$,
and expanding this
out as in the proof of Proposition~\ref{Lp-geod} we find
\[
1 + |\gamma|^2 - 2 \mbox{Re}\, (\bar{\gamma}\zeta) =
\sum_{j=1}^\infty |c_j|^p (1 + |\alpha_j|^2 - 2 \mbox{Re}\,
(\bar{\alpha}\zeta) ) \qquad (|\zeta| = 1).
\]
As both sides are harmonic for $\zeta \in \Bbb D$ and continuous on
$\bar{\Bbb D}$, the equality remains valid for $\zeta \in \Bbb D$. Putting
$\zeta=0$ gives $\sum_j |c_j|^p (1 + |\alpha_j|^2) = 1 +|\gamma|^2$ and
then it follows that $\sum_j |c_j|^p \alpha_j = \gamma$. \qed

\begin{lemma}
\label{gent-lem}
Let $f \in H^\infty({\Bbb D})$, $\gamma \in {\Bbb D}$. Then

\[
\displaystyle \frac{f(e^{i\theta})}
{\left(\displaystyle
\frac{e^{i\theta}-\gamma}{1- \bar{\gamma}e^{i\theta}}\right)}
\]
is a
non-negative real number for almost all $\theta \in \IR$ if and only if
\[
f(\zeta) = t
\left( \frac{\zeta - \alpha}{1 - \bar{\alpha} \zeta} \right)
\left( \frac{1 - \bar{\alpha} \zeta}{1- \bar{\gamma}\zeta} \right)^2
\]
for some $t \geq 0$, $|\alpha| \leq 1$.
\end{lemma}

\noindent
{\sc Proof}: The case $\gamma=0$ is due to Gentili
\cite{gent-l1} (see lemma 2 and the
proof of theorem 6). The general case follows from Gentili's result
by the change of variables
$\eta = (\zeta-\gamma)/(1-\bar{\gamma}\zeta)$. \qed

\begin{lemma}
\label{l1-uniq}
If $X= \ell^1_n$ or $X = \ell^1$
and $x \neq y \in B_X$, then there is a unique normalized
complex geodesic in $B_X$ joining $x$ and $y$.
\end{lemma}

\noindent
{\sc Proof}: As already noted,
we know that there exists a normalized complex geodesic
$\phi$ in $B_X$ joining $x$ and $y$. In fact we have an explicit form
(\ref{geod-form-N}) of one such $\phi = (\phi_j)_j$ by
Proposition~\ref{lp-exist}.
, where
\[
\phi_j(\zeta) = c_j
\left( \frac{\zeta - \alpha_j}{ 1 - \bar{\alpha}_j \zeta \zeta}
\right)^{\beta_j}
\left( \frac{ 1 - \bar{\alpha}_j \zeta}{ 1 - \bar{\gamma} \zeta} \right)^2
\]
for all $j$

Now suppose $\psi$ is
a second normalized complex geodesic joining $x$ and $y$.
Suppose $\phi(0)= \psi(0) = x$ and $\phi(s) = \psi(s) = y$ where $s >0$.
The argument given earlier in the proof of Theorem~\ref{uniq-arnp-th}
shows that $f = \lambda\phi + (1-\lambda)\psi$ is also a complex geodesic
joining $x$ and $y$.
Thus $\|f(e^{i\theta})\| = 1$ for all $\theta \in \IR$. It follows
that we must have equality in the triangle inequality
\[
1 = \|f(e^{i\theta})\| = \sum_j | \lambda \phi_j(e^{i\theta})
+(1-\lambda) \psi_j(e^{i\theta})| \leq \sum_j
 \lambda| \phi_j(e^{i\theta})| + (1-\lambda) |\psi_j(e^{i\theta})| = 1.
\]
This forces
\[
\frac{ \psi_j(e^{i\theta}) }{\phi_j(e^{i\theta})}
\]
to be a non-negative real number for all $\theta$ except those for
which the denominator is zero.

To prove that $\psi_j = \phi_j$, we consider the case $c_j =0$ and $c_j \neq
0$ separately.
In the first case $\phi_j \equiv 0$ and hence
$x_j = y_j = 0$. Since linear isometries of $X$ map complex geodesics in
$B_X$ to complex geodesics, $\psi$ and
\[
\tilde{\psi} = (\psi_1, \ldots, \psi_{j-1}, - \psi_j, \psi_{j+1}, \ldots)
\]
are both normalized complex geodesics joining $x$ and $y$. Therefore, so is $g =
(\psi + \tilde{\psi})/2$. Since $\|g(e^{i\theta})\| =
\|\psi(e^{i\theta})\| = 1$ for all $\theta$, it easily follows that $\psi_j
\equiv 0 \equiv \phi_j$.

If $c_j \neq 0$, 
Lemma~\ref{gent-lem} applied to the function
\[
\frac{ \phi_j(\zeta) }{ c_j }
\left( \frac{ \zeta - \alpha_j}{1 - \bar{\alpha}_j \zeta } \right)^{1 -
\beta_j}
\]
shows that
\[
\left.
\frac{ \phi_j(e^{i\theta}) }{ c_j }
\left( \frac{ e^{i\theta} - \alpha_j}{1 - \bar{\alpha}_j e^{i\theta} }
\right)^{1 - \beta_j}
\right/
\left( \frac{ e^{i \theta} - \gamma}{ 1 - \bar{\gamma} e^{i\theta} }\right)
\]
is non-negative for almost all $\theta$. Therefore the same holds for
$\psi_j(e^{i\theta})$ in place of $\phi_j(e^{i\theta})$
and Lemma~\ref{gent-lem} then shows that
there exist $t > 0$ and $|\beta| \leq 1$ so that
\[
\psi_j(\zeta) = t c_j
\left( \frac{ \zeta - \alpha_j}{1 - \bar{\alpha}_j \zeta }
\right)^{\beta_j-1}
\left(\frac{ \zeta - \beta}{1 - \bar{\beta}\zeta} \right)
\left( \frac{1 - \bar{\beta}\zeta}{1 - \bar{\gamma}\zeta} \right)^2.
\]

If $\beta_j = 0$ and $|\alpha_j| < 1$, then analyticity of $\psi_j(\zeta)$
forces $\alpha_j = \beta_j$ and two terms cancel in
the expression for $\psi_j(\zeta)$. Then from $\phi_j(0) = \psi_j(0) =
x_j$, we conclude that $c_j = t c_j$ and $t = 1$. So $\phi_j = \psi_j$ in
this situation.

In the remaining cases, $\phi_j(0) = \psi_j(0) = x_j$
yields
\[
t c_j (-\alpha_j)^{\beta_j -1} (-\beta) + c_j (- \alpha_j)^{\beta_j}
\]
and hence $\alpha_j = t \beta$.
Then $\psi_j(s) = \phi_j(s) = y_j$ and some cancellation of common terms
shows that
\[
(s - \alpha_j)(1 - \bar{\alpha}_j s)  = t(s - \beta)(1 - \bar{\beta} s).
\]
Combining this with $\alpha_j = t \beta$ yields
\[
(1-t)s(1 - t |\beta|^2) = 0.
\]
Hence $t = 1$ or $t = |\beta|^{-2}$.
In the second case,
\[
\alpha_j = t \beta = 1/\bar{\beta}.
\]
Since $|\alpha_j| \leq 1$ and $|\beta| \leq 1$, we must
have $|\alpha_j| = |\beta| = 1$ so that $t=1$. Hence
$\alpha_j = \beta$ and $\psi_j (\zeta) = \phi_j (\zeta)$.
\qed

We can summarise our results for $\ell^p$ as follows.

\begin{corollary}
Let $B_p$ denote the unit ball of $\ell^p$, $1 \leq p < \infty$. Then
\begin{description}
\item{(i)} Any two distinct points in $B_p$ can be joined by a unique
normalized complex geodesic.

\item{(ii)} All complex geodesics in $B_p$ are continuous.

\item{(iii)} A map $\phi\colon {\Bbb D} \to B_p$ is a complex geodesic
if and only if it is a non-constant map of the form given in
Proposition~\ref{lp-exist}
\end{description}
\end{corollary}

\noindent
{\sc Proof}: For all $1 \leq p < \infty$, existence follows from
Theorem~\ref{1-compl} (or from \cite{DTV}). Uniqueness for $p > 1$ follows
from Theorem~\ref{uniq-arnp-th} and uniqueness for $p=1$ has just been
established in Lemma~\ref{l1-uniq}.

(ii) follows from Theorem~\ref{th-cont}. For $p=1$, this has already
been noted in Examples~\ref{uniq-ex}. It is straightforward
that Theorem~\ref{th-cont} applies to $\ell^p$ for $1 < p < \infty$
because $\ell^p$ is uniformly convex in the real sense (see \cite{DGT}).

(iii) follows from Proposition~\ref{lp-exist} and (i). \qed

We suspect that a more general version of this result holds for $L^p(\mu)$
in place of $\ell^p$, but we have not managed to prove $L^p$-versions of
Propositions \ref{lp-exist} or \ref{l1-uniq}.

\begin{example}
\label{lp+lq}
Let $X = \ell^{p_1} \oplus_r \ell^{p_2} = \{ x = (y,z):
y \in \ell^{p_1}, z \in \ell^{p_2} \}$
normed by $\|x\| = \left( \|y\|_{p_1}^r + \|z\|_{p_2}^r \right)^{1/r}$.

One can check using Proposition~\ref{lempert-criterion} that
for $1 \leq p_i < \infty$, $1 \leq r < \infty$ all nonconstant maps
$\phi \colon {\Bbb D} \to B_X$ of the following form are complex geodesics.
\begin{eqnarray*}
\phi(\zeta) &=& ( (\phi_{1j})_{j=1}^\infty, (\phi_{2j})_{j=1}^\infty )\\
\phi_{ij}(\zeta) &=& c_{ij} \left(
\frac{\zeta - \alpha_{ij}}{ 1 - \bar{\alpha}_{ij} \zeta}\right)^{\beta_{ij}}
\left(
\frac{ 1 - \bar{\alpha}_{ij} \zeta}{ 1 - \bar{\gamma}_i \zeta}
\right)^{2/p_i}
\left( \frac{ 1 - \bar{\gamma}_i \zeta}
{1- \bar{\gamma}\zeta} \right)^{2/r}
\end{eqnarray*}
where $|\alpha_{ij}| \leq 1$, $|\gamma_i| < 1$, $|\gamma| < 1$,
$\beta_{ij}$ is $0$ or $1$, and the following relations hold
\begin{eqnarray*}
\sum_j |c_{ij}|^{p_i} \alpha_{ij} &=& \gamma c_i^{p_i} \qquad (i = 1,2)\\
c_1^r \gamma_1 + c_2^r \gamma_2 &=& \gamma\\
c_1^r (1 + |\gamma_1|^2) + c_2^r (1 + |\gamma_2|^2) &=& 1 + |\gamma|^2
\end{eqnarray*}
where
\[
c_i = \left( \frac{1}{1+|\gamma_i|^2} \sum_j |c_{ij}|^{p_i}
(1 + |\alpha_{ij}|^2) \right)^{1/p_i}
\quad (i = 1,2).
\]

The proof of this involves observing first that for
$x =(y,z) \in X$ with $\|x\|=1$,
\[
N_x = ( \|y\|^{r-1} N_{y/\|y\|}, \|z\|^{r-1} N_{z/\|z\|} )
\]
with $N_{y/\|y\|}$ and $N_{z/\|z\|}$ given as in the proof of
Proposition~\ref{Lp-geod} for $\ell^p$.
To apply Proposition~\ref{lempert-criterion}, take $p(\zeta) =
|1 - \bar{\gamma}\zeta|^2$ and
\[
h(\zeta) = \left( \tilde{c}_{ij} (1 - |\alpha_{ij}|^2 )^{2-2/p_i}
(1 -\bar{\gamma_i} \zeta)^{2/p_i - 2/r} (1 - \bar{\gamma}\zeta)^{2/r}
\left( \frac{\zeta - \alpha_{ij} }{1 - \bar{\alpha}_{ij} \zeta }
\right)^{1 - \beta_{ij} } \right)_{ij}
\]
where
\[
\tilde{c}_{ij} = c_i^{r - p_i} |c_{ij}|^{p_i - 2} \bar{c}_{ij}
\quad(i=1,2; 1 \leq j < \infty).
\]

We suspect that all complex geodesics in $B_X$ are of this form.
Other examples of complex geodesics in spaces which are direct sums
of more than two summands of $\ell^p$-type can also be exhibited.

\end{example}
\begin{remark}
\label{l-infty}
The case $p = \infty$ is excluded in all of the previous calculations
because it is well known that
almost everything is different for $\ell^\infty$ (and
$L^\infty$). Even for the unit ball $B_{\infty,2} = \{ (z_1,z_2) \in
\IC^2 : \max_i |z_i| < 1 \}$ (polydisc) of $\ell^\infty_2$, many of
the differences are apparent. The only points of $\partial
B_{\infty,2}$ that are complex extreme points are those where $|z_1|
= |z_2| = 1$ and therefore the result of Vesentini~\cite{vesentini}
cited at the beginning of Section~\ref{sec-uniq} shows that there are
many complex geodesics joining $0$ to $z = (z_1,z_2)$ if $|z_1| \neq
|z_2|$. In fact, if $|z_2| < |z_1|$, the normalized complex geodesics
joining $0$ to $(z_1,z_2)$ are
\[
\phi(\zeta) = \left( \zeta \frac{z_1}{|z_1|}, g(\zeta) \right)
\]
where $g$ is any analytic function on $\Bbb D$ with $g(0) = 0$,
$g(|z_1|) = z_2$ and $\sup_{\zeta \in \Bbb D} |g(\zeta)| < 1$. (This
can easily be verified using the fact that the Kobayashi distance
between $x,y \in B_{\infty,2}$ is given by
$\max_{i=1,2} \rho (x_i, y_i)$.) Thus we see non-uniqueness and
discontinuity of complex geodesics. Since $B_{\infty,2}$ is a
homogeneous domain, for any pair of points $x,y \in B_{\infty,2}$, we
can find a biholomorphic automorphism $F$ with $F(0) = x$ (in fact 
$F(z_1,z_2) = ( (x_1 - z_1)/(1- \bar{x}_1 z_1), (x_2 - z_2)/(1 -
\bar{x}_2 z_2) )$ will do). Then for any normalized complex geodesic
$\phi$ joining $0$ and $F^{-1}(y)$, $F \circ \phi$ is a normalized
complex geodesic joining $x$ and $y$.

The infinite dimensional case $\ell^\infty$ is somewhat similar to
the finite dimensional one. We state without proof the following
result ($B_\infty$ is the unit ball of $\ell^\infty$).
\begin{quote}
{
A map $\phi = (\phi_1, \phi_2, \ldots ) \colon {\Bbb D} \to B_\infty$
is holomorphic if and only if each coordinate function
$\phi_j(\zeta)$ is analytic.}
\end{quote}
One can check that $\phi$ is a complex geodesic if and only if either
there exists $j$ for which $\phi_j$ is an automorphism of $\Bbb D$ or
else there is a subsequence of $(\phi_j)_j$ which converges to an
automorphism. This can be done directly or via a result of
Gentili~\cite{gent-bal} who gives a necessary and sufficient conditions
for a holomorphic map $\phi \colon X \to B_X$ to be a complex
geodesic in the case when $X = C(K)$ is the space of continuous
functions on a compact Hausdorff space $K$. The condition is that
there exists $k \in K$ so that $\zeta \mapsto \phi(\zeta)(k)$ is an
automorphism of the unit disc. Since $\ell^\infty$ is the same as the
space of continuous functions on the Stone-\v{C}ech compactification
of the integers, we can apply his result.
 There will be more than one complex geodesic joining
$x,y \in B_\infty$ unless $|(x_j - y_j)/(1 - \bar{x}_j y_j)|$ is
constant.

Similar remarks apply to $c_0$. In this case there are no complex
extreme points on the unit sphere and there is more than one
normalized complex
geodesic joining every pair of points in the unit ball.
\end{remark}
\section{TAUTNESS AND CURVATURE}
\setcounter{equation}{0}
\label{sec-taut}

In previous sections, we have used (for convex domains) convergence
principles to establish the existence of holomorphic mappings with certain
extremal properties. These convergence properties have been formalised for
finite-dimensional domains and manifolds and certain relationships
established. In this section, we show that these results do not extend to
arbitrary Banach spaces, even for convex bounded domains and give a result
on curvature of the Kobayashi infinitesimal metric.

We define hyperbolic, complete hyperbolic and taut domains in a Banach
space. The first two definitions are standard, while various versions of
the third are possible.

\begin{definition}
\label{def-hyperbolic}
A domain $\cal D$ in a Banach space $X$ is called {\em hyperbolic} if
$K_{\cal D}$ induces the original topology on $\cal D$. If, moreover,
$({\cal D}, K_{\cal D})$ is a complete metric space, we call $\cal D$ a
{\em complete hyperbolic domain}.
\end{definition}

If $\cal D$ is hyperbolic, then $K_{\cal D}$ is a distance on $\cal D$
({\em i.e.} $K_{\cal D}$ separates the points of $\cal D$) and
the converse is true for finite dimensional domains (see
\cite{royden,harris}). We do not know if  the converse is true for
infinite dimensional Banach spaces. Harris~\cite[theorem 24]{harris}
has proved that
a convex bounded domain in a Banach space is complete hyperbolic
and Barth~\cite{barth} has proved that a convex  domain in $\IC^n$
which contains no complex lines is hyperbolic (and indeed
biholomorphically equivalent to a bounded domain). Again, we do not
know if this result extends to infinite dimensional Banach spaces.

\begin{definition}
\label{def-taut}
A domain $\cal D$ in a Banach space $X$ is {\em taut} if there exists
a Hausdorff locally convex topology $\tau$ on $X$
such  that every net $(f_\alpha)_\alpha$ in
$H({\Bbb D}, {\cal D})$ contains either a compactly divergent
subnet or a subnet which converges with respect to $\tau$, uniformly
on compact subsets of $\Bbb D$, to some $f \in H({\Bbb D}, {\cal D})$.

A net $(f_\alpha)_{\alpha \in \Gamma}$ is compactly divergent
if, given $K$ compact in $\Bbb D$ and $L$ a norm-compact subset
of $\cal D$, there exists $\alpha_0 \in \Gamma$ such that
$f_\alpha(K) \cap L = \emptyset$ for all $\alpha \geq \alpha_0$.
\end{definition}

For $\cal D$ finite dimensional, $\tau$ has to be the same as the norm
topology and we have
\[
\mbox{complete hyperbolic} \Rightarrow \mbox{taut} \Rightarrow
\mbox{hyperbolic}
\]
(and the converses are both false \cite{kiernan,barth-counter}).
Also, for finite dimensional domains
the unit disc $\Bbb D$ can be replaced by any finite dimensional
domain in the definition of tautness \cite{barth-taut,wu}.

\begin{proposition}
\label{taut-dual}
The unit ball $B_X$ of a Banach space $X$ is taut if and only if $X$
is isometrically isomorphic to a dual space.
 
Moreover, if
$X$ contains any bounded convex domain $\cal D$ which is taut
(for a locally convex topology $\tau$) then $X$ is
isomorphic to a dual Banach space and
the topology $\tau$ is weaker than the norm topology.
\end{proposition}

\noindent
{\sc Proof}: We first suppose that $B_X$ is taut and that
$\tau$ is a locally convex Hausdorff
topology associated with tautness.
Let $(x_\alpha)_\alpha$ be a net in $\bar{B}_X$ and for each $\alpha$
let $f_\alpha(\zeta) = \zeta x_\alpha $ for $\zeta \in \Bbb D$. Then
$f_\alpha$ is a net in $H({\Bbb D}, B_X)$ and, since $f_\alpha(0) =
0$ for all $\alpha$, it contains no compactly divergent subnet.
Therefore it  must have a subnet $(f_\beta)_\beta$ which converges
relative to $\tau$, uniformly on compact subsets of $\Bbb D$, to some
$f \in H({\Bbb D}, B_X)$. Since $\tau$ is a locally convex topology,
it is easily seen that $f(\zeta) = \zeta x$ for some $x \in X$.
Because $\|f(\zeta)\| < 1$ for all $\zeta \in \Bbb D$, we must have
$\|x\| \leq 1$ and moreover $x_\beta \to x$ (with respect to $\tau$).
This shows that $(\bar{B}_X, \tau)$ is compact.
A result of Ng~\cite{ng} now implies that $X$
is isometrically isomorphic to a dual space.

Conversely, suppose now that $X$ is isometrically isomorphic to the
dual of a Banach space $Y$. Then $\bar{B}_X$ is compact in the weak*-
or $\sigma(X,Y)$-topology. Let $(f_\alpha)_\alpha$ denote a net in
$H({\Bbb D}, B_X)$. Consider the Taylor series expansions
\[
f_\alpha(\zeta) = \sum_n a_{\alpha,n} \zeta^n.
\]
From the Cauchy formula, we have $\|a_{\alpha,n}\| \leq 1$  for all
$\alpha$ and $n$. By compactness of the product of infinitely many
copies of $\bar{B}_X$ (for the product $\sigma(X,Y)$-topology) we can
find a subnet $(f_\beta)_\beta$ so that $a_{\beta,n} \to a_n \in
\bar{B}_X$ for all $n$.

If $\|a_0\| = 1$, then for any $0 < r < 1$,
$f_\beta(0) = a_{\beta,n} \not \in r \bar{B}_X$ if $\beta$ is
sufficiently large. Thus $\|f_\beta(0)\| \to 1$. Now
\begin{eqnarray*}
\tanh^{-1} \| f_\beta(\zeta) \| &=& K_X (f_\beta(\zeta),0)\\
&\geq& K_X( f_\beta(0), 0 ) - K_X(f_\beta(\zeta), f_\beta(0))\\
&\geq& \tanh^{-1} \|f_\beta(0) \| - \rho(\zeta,0)
\end{eqnarray*}
(here we have used
the distance decreasing property for the Kobayashi metric)
shows that $\| f_\beta(\zeta) \| \to 1$
uniformly on compact subsets of $\Bbb D$. Thus $f_\beta$ is
compactly divergent.

For the case $\|a_0\| < 1$, let $f(\zeta) = \sum_n a_n \zeta^n$.
Since $\|a_n\| \leq 1$ for all $n$, clearly $f \in H({\Bbb D}, X)$.
Fix $y \in Y$, $0 < r < 1$ and $\epsilon > 0$. Choose $N$ so that
$\sum_{n=N+1}^\infty r^n < \epsilon$. Then, by considering
power series, we obtain
\begin{eqnarray*}
\sup_{|\zeta| \leq r} | \langle y, f_\beta(\zeta) - f(\zeta) \rangle|
&\leq& \sum_{n=0}^N |\langle y, a_{\beta,n} - a_n \rangle| +
2\sum_{n=N+1}^\infty r^n \|y\| \\
&\leq& \sum_{n=0}^N | \langle y, a_{\beta,n} - a_n \rangle| +
2\epsilon \|y\|
\end{eqnarray*}
for all $\beta$. Hence $f_\beta \to f$ with respect to $\sigma(Y,X)$,
uniformly on compact subsets of $\Bbb D$.

Since $f_\beta(\zeta) \in B_X$ for all $\zeta$ and $\beta$, it
follows that $f$ has values in $\bar{B}_X$. Since $f(0) = a_0 \in
B_X$, $f \in H({\Bbb D}, B_X)$. This completes the proof that $B_X$
is taut when $X$ is a dual space.

It remains to show the last part of the proposition. Suppose ${\cal D}$
is a taut bounded convex domain in $X$.
We can assume that $0 \in \cal D$
by a translation. A simple argument based on compactness of the unit circle
shows that $ {\cal D}_0 =
\cap_{\theta \in \IR} e^{i\theta} \cal D$ is a bounded balanced
convex domain. By replacing the norm on $X$ with an equivalent one, we can
assume that ${\cal D}_0 = B_X$. The same argument as used at the
beginning of the proof, together with the observation that a map $f
\colon {\Bbb D} \to X$ of the form $f(\zeta) = \zeta x$ has its
values in $\cal D$ if and only if it has values in ${\cal D}_0 =
B_X$, shows that $X$ is a dual space. Moreover the proof of
Ng~\cite{ng} implies that $\tau$ is weaker than the norm topology.
He shows that $X$ is the dual of a certain subspace
$Y \subset X^*$ and that the identity map
$id \colon (\bar{B}_X, \tau) \to (\bar{B}_X, \sigma (X,Y) )$ is
a homeomorphism. It follows that $id \colon (\bar{B}_X , \| \cdot \|)
\to (\bar{B}_X, \tau)$ is continuous, and hence that $id \colon (X,
\| \cdot \|) \to (X, \tau)$ is continuous.
\qed

We remark that we could have appealed to this proposition for the
proof of Proposition~\ref{lp-exist} in place of the more direct
argument we used. 

\begin{example}
For $X = c_0$, $B_X$ is complete hyperbolic (by \cite{harris})
but not taut.
\end{example}

The second part of the proof of Proposition~\ref{taut-dual} can
easily be modified to show that the following are taut
\begin{description}
\item{(a)} convex bounded domains in reflexive Banach spaces (for
$\tau$ the weak topology);

\item{(b)} convex bounded domains in dual Banach spaces with the
property that their norm closures are weak*-compact ($\tau$ the
weak*-topology).
\end{description}

\begin{proposition}
\label{taut-con}
If $\cal D$ is a convex bounded taut domain in a Banach space $X$, then
$\cal D$ is C-connected.
\end{proposition}

\noindent
{\sc Proof}: Given two points
$x,y \in \cal D$, choose $(f_n)_n$ in $H({\Bbb D}, {\cal
D})$ so that $f_n(0) = x$, $f_n(s_n) = y$, $s_n \in (0,1)$ for all $n$ and
$s_n \to s = \tanh K_{\cal D} (x,y)$ as $n \to \infty$.

Since $\cal D$ is taut and $f_n(0) = x$ for all $n$, it follows that
$(f_n)_n$ has a $\tau$-convergent subnet ($\tau$ being the locally convex
topology related to tautness of $\cal D$). Let $f \in H({\Bbb D}, {\cal
D})$ denote the limit of one such subnet. Clearly $f(0) = x$ and we claim
$f(s) = y$ (which will show that $f$ is a complex geodesic by convexity and
Proposition~\ref{geod-equiv}).

For $r = (1+s)/2$, we have
\[
f_n(s_n) - f_n(s) = \frac{s_n -s}{2\pi i} \int_{|z|=r}
\frac{f_n(z)}{(z-s_n)(z-s)} \, dz.
\]
Since $\cal D$ is bounded and $s_n \to s$, it follows that $f_n(s_n)
-f_n(s) \to 0$ in norm (hence also in the topology $\tau$ by
Proposition~\ref{taut-dual}) as $n \to
\infty$. Now if $(f_{n_\alpha})_{\alpha \in \Gamma}$ is a subnet of
$(f_n)_n$ which converges to $f$ uniformly with respect to $\tau$ on compact
subsets of $\Bbb D$, then
\[
f_{n_\alpha} (s_{n_\alpha}) - f(s) = 
f_{n_\alpha} (s_{n_\alpha}) - f_{n_\alpha} (s) + f_{n_\alpha} (s) - f(s) 
\to 0
\]
in the topology $\tau$. Since $f_{n_\alpha} (s_{n_\alpha}) = y$, it follows
that $f(s) = y$. \qed

We observe that one can similarly prove an infinitesimal version of
Proposition~\ref{taut-con}:
{\em if $\cal D$ is a convex bounded taut domain in a Banach space $X$,
$x \in \cal D$ and $0 \neq v \in X$, then there exists a complex
geodesic $f \in H({\Bbb D}, {\cal D})$ such that $f(0) = x$ and
$f'(0) = v/k_{\cal D}(x,v)$.}

\begin{definition}
\label{curv-def}
If $\cal D$ is a bounded domain in a Banach space $X$, then the
{\em holomorphic sectional curvature} of the Kobayashi infinitesimal metric
$k_{\cal D}$ at a point $x \in \cal D$ in the (nonzero)
direction $v \in X$ is
\[
\kappa_k(x,v) = \sup \left\{
\frac{\triangle \log k_{\cal D}^2(f(z), f'(z))\mid_{z=0}}
{- 2 k_{\cal D}^2 (f(0), f'(0))} :
f \in H(r{\Bbb D}, {\cal D}), r > 0, f(0) = x, f'(0) = v \right\}.
\]

Since $k_{\cal D}$ may not even be continuous in general, $\triangle$ above
denotes the generalised Laplacian. This is defined (motivated by
\cite{heins}, see \cite{burbea}) for upper semicontinuous
functions $u$ with values in $[-\infty, \infty)$ (but not identically
$-\infty$) at points $z$ with $u(z) \neq -\infty$ by
\[
\triangle u(z) = 4 \liminf_{r \to 0^+} \frac{1}{r^2} \left\{
\frac{1}{2\pi} \int_0^{2\pi} u(z + r e^{i\theta}) \, d\theta - u(z)
\right\}.
\]
\end{definition}

At a local maximum point $z$, $\triangle u (z) \leq 0$ and at a finite
local minimum $\triangle u (z) \geq 0$.

The following result for finite-dimensional domains is due to
Wong~\cite{wong}, Burbea~\cite{burbea} and
Suzuki~\cite{suz-hol,suz-int}.

\begin{proposition}
\label{curv-4}
If $\cal D$ is a taut bounded convex domain in a Banach space $X$, then
\[
\kappa_k (x,v) = -4 \qquad (x \in {\cal D}, 0 \neq v \in X).
\]
\end{proposition}

\noindent
{\sc Proof}: We will use the fact that this result is true for the unit
disc $\Bbb D$ (see for instance \cite[\S 3.4]{dineen}) --- in fact every
suitable function $f$ attains the supremum in Definition~\ref{curv-def}
when ${\cal D} = \Bbb D$ and the supremum is $-4$.

For general $\cal D$ and
$f$ as in the definition of $\kappa_k(x,v)$, let
$\lambda_f(z) = c_{\cal D}^2 (f(z), f'(z))$ (which coincides with 
$k_{\cal D}^2 (f(z), f'(z))$ by convexity of $\cal D$). By Montel's
theorem, we can find $g \in H({\cal D}, {\Bbb D})$ such that $g(x) = 0$ and
\[
c_{\cal D} (x,v) = |g'(x)(v)| = c_{\Bbb D} (g(x), g'(x)(v) ) =
\alpha(g(x), g'(x)(v)).
\]

Let $\beta(z) = \alpha^2( g(f(z)), g'(f(z))(f'(z)) )$ for $z \in r \Bbb
D$. Observe that
\[
\beta(z) \leq c_{\cal D}^2 (f(z), f'(z)) = \lambda_f(z)
\]
and $\beta(0) = \alpha^2( g(x), g'(x)(v)) = c_{\cal D}^2 (x,v) =
\lambda_f(0)$. Since $v \neq 0$ and $\cal D$ is bounded, there is a
neighbourhood of $0$ where $\beta(z)$ does not vanish. Hence $\log
(\lambda_f/\beta )$ has a local minimum at the origin and therefore
\[
\left. \triangle \log \frac{\lambda_f}{\beta} \right|_{z=0}
= \triangle \log \lambda_f (0) - \triangle \log \beta(0)  \geq 0,
\]
where we have used the fact that $\log \beta$ is twice continuously
differentiable (and hence the $\liminf$ in the definition of
$\triangle \log \beta$ is a limit).
It follows that
\[
\frac{ \triangle \log \lambda_f (0) }{-2\lambda_f(0)} \leq
\frac{\triangle \log \beta (0)}{-2 \beta(0)} = -4
\]
and hence that $\kappa_k (x,v) \leq -4$.

Since $\cal D$ is taut and bounded, there exists a complex geodesic
$f \in H({\Bbb D}, {\cal D})$ such that $f(0) = x$ and $f'(0) = v/k_{\cal
D}(x,v)$ (by the infinitesimal version of Proposition~\ref{taut-con}). For
this $f$ we have
\[
k_{\cal D}^2 ( f(z), f'(z) ) = \alpha^2 (z, 1) = \frac{1}{1 - |z|^2}
\]
(see Proposition~\ref{geod-equiv})
and the one variable result shows that $f$ attains the upper bound of $-4$.
\qed

\begin{corollary}
If $X$ is a dual Banach space, then the Kobayashi infinitesimal metric $k_X$
on $B_X$ has constant holomorphic sectional curvature $-4$.
\end{corollary}

\vfill

This paper appeared in Progress in Functional Analysis, K.D. Bierstedt, J.
Bonet, J. Horv\'ath \& M. Maestre (Eds.) pp 333--365, Elsevier (1992) .

\vfill

\end{document}